\PassOptionsToPackage{capitalize, sort, compress}{cleveref}
\documentclass[11pt,a4paper]{article}
\usepackage{color}
\usepackage[colorlinks=true]{hyperref}
\usepackage[margin=1in]{geometry}
\usepackage{fullpage,mathrsfs,graphicx,framed,color,amssymb,amsmath,amsthm}
\usepackage[boxed, noline, ruled, linesnumbered]{algorithm2e}
\usepackage{epsfig, graphicx}
\usepackage{latexsym,amsfonts,amsbsy,amssymb}

\usepackage{lmodern}

\usepackage{amsmath,amsthm}
\usepackage{enumerate}
\usepackage{cleveref}
\usepackage{tikz}
\usepackage{stmaryrd} 
\usepackage{pgfplots}
\pgfplotsset{compat=1.17}
\usepackage{multirow}
\usepackage{booktabs}
\usepackage{arydshln} 
\usepackage{bbding} 

\crefname{equation}{}{}
\crefname{theorem}{Theorem}{Theorems}
\crefname{figure}{Figure}{Figures}
\crefname{table}{Table}{Tables}
\crefname{section}{Section}{Sections}
\crefname{example}{Example}{Examples}
\crefname{lemma}{Lemma}{Lemmas}
\crefname{remark}{Remark}{Remarks}

\usepackage{subcaption}
\usepackage{amsmath}    
\usepackage{amssymb}    
\usepackage{changes} 
\usepackage{tikz}
\usetikzlibrary{arrows.meta,positioning,calc}
\usepackage[utf8]{inputenc} 

\usepackage{caption}
\usepackage{graphicx} 
\usepackage[capitalize, sort, compress]{cleveref}
\usepackage{hyperref} 

\makeatletter 

\@addtoreset{equation}{section}
\makeatother 
\allowdisplaybreaks 

\newtheorem{remark}{Remark}[section]

\allowdisplaybreaks 

\numberwithin{equation}{section}

\begin{document}
\title{Quadrature-Enhanced Monte Carlo fPINN Method for High-Dimensional Fractional PDEs}
	
\author{Qingkui Ma\footnote{School of Mathematics and Statistics \& 
			Hubei Key Laboratory of Mathematical Sciences,
			Central China Normal University, Wuhan 430079, China
			(maqingkui@mails.ccnu.edu.cn).}\ \ \ 
		Hehu Xie\footnote{SKLMS, NCMIS, Academy of Mathematics and Systems Science,
			Chinese Academy of Sciences, No.55, Zhongguancun Donglu, Beijing 100190, 
			China, and School of Mathematical Sciences, University of Chinese Academy
			of Sciences, Beijing, 100049, China (hhxie@lsec.cc.ac.cn).}\ \ \ and \ \ 
		Xiaobo Yin\footnote{School of Mathematics and Statistics \& Key Laboratory of Nonlinear Analysis \& Applications (Ministry of Education),
			Central China Normal University, Wuhan 430079, China (yinxb@ccnu.edu.cn).}}
\date{}
\maketitle
	
\begin{abstract}
Fractional PDEs involving the fractional Laplacian on bounded domains are challenging because of hypersingular nonlocal kernels, exterior Dirichlet constraints, reduced boundary regularity, and the high computational cost in high dimensions. To address these issues, we first adopt a spatially varying radius with directional distance-to-boundary information, which yields a geometry-adaptive three-part decomposition of the fractional Laplacian: singular near-field, regular interior far-field, and analytical exterior far-field contributions.
Then we employ Gauss-Jacobi quadrature for the singular radial integral, Gauss quadrature for the regular interior radial integral, and Monte Carlo sampling for the angular variables. A feature-enhanced physics-informed neural network trial space is finally used to tackle the low-regularity behavior near the boundary. Through the above steps, we obtain a quadrature-enhanced Monte Carlo fractional physics-informed neural network (QE-MC-fPINN) method. Numerical experiments on fractional Poisson equations and time-dependent fractional PDEs show that, on the tested benchmarks, the proposed method outperforms two representative MC-fPINN discretizations in accuracy and convergence, especially for solutions with strong boundary singularities. 

\end{abstract}
	
\section{Introduction}
Fractional partial differential equations (PDEs) involving the fractional Laplacian arise in anomalous diffusion, long-range interaction, and nonlocal transport models \cite{Bonito2018,DElia2020,Lischke2018,Caffarelli2007}. They also appear in several high-dimensional settings, including kinetic and fractional Fokker--Planck models posed in phase space, mean field games driven by nonlocal L\'evy-type diffusions, and multi-asset option pricing under jump processes \cite{guo2022monte,ersland2021fractional,guo2018tempered}. In such applications, the dimension often arises from phase-space variables, coupled state variables, or parametric uncertainty.
The fractional Laplacian admits several distinct, generally non-equivalent definitions depending on the underlying setting \cite{Lischke2018}. In this paper, we adopt the Riesz definition with zero exterior extension:
\begin{equation}\label{eq:fractional_laplacian}
		(-\Delta)^{\alpha/2} u(\boldsymbol{x})
		=
		C_{d,\alpha}\,\mathrm{P.V.}\!\int_{\mathbb{R}^d}
		\frac{u(\boldsymbol{x})-u(\boldsymbol{y})}
		{\|\boldsymbol{y}-\boldsymbol{x}\|^{d+\alpha}}
		\,d\boldsymbol{y},
		\qquad
		C_{d,\alpha}
		:=
		\frac{2^{\alpha-1}\alpha\,\Gamma\!\left(\alpha/2+d/2\right)}
		{\pi^{d/2}\Gamma(1-\alpha/2)},
\end{equation}
where $\mathrm{P.V.}$ denotes the Cauchy principal value and $0<\alpha<2$. Throughout this work, we focus on bounded-domain problems with exterior Dirichlet conditions. This setting is of particular interest because the interaction between nonlocality and the boundary strongly affects both the regularity of the solution and the design of numerical methods. Such bounded-domain fractional models arise in anomalous transport, subsurface flow, and nonlocal diffusion processes in engineering systems \cite{Lischke2018,pang2019fpinns,hu2024tackling}. The numerical treatment of these problems is considerably more difficult than that of classical local PDEs. The main difficulties stem from the hypersingular and nonlocal nature of the operator, which induces global coupling, and from the reduced regularity of solutions near the boundary on bounded domains \cite{Lischke2018}. These challenges are further amplified when exterior Dirichlet constraints must be imposed accurately and when the problem dimension is high.
	
Classical discretizations include finite difference and finite element methods \cite{Huang2023AGF,Meerschaert2006FDM,Acosta2017,Ainsworth2018,Gao2019AFE,Sheng2024}. While effective in low-dimensional settings, such methods may require fine meshes and expensive nonlocal matrix operations, especially for high-dimensional problems or domains with complex boundary geometry. In recent years, neural-network-based solvers have provided a complementary mesh-free framework. Starting from early neural-network methods for differential equations \cite{lagaris1998artificial}, physics-informed neural networks (PINNs) \cite{raissi2019physics} and their fractional variants, including fPINNs \cite{pang2019fpinns}, bi-orthogonal fPINNs \cite{ma2023bi}, and spectral-fPINNs \cite{zhang2025spectral}, have shown promising performance for forward and inverse fractional PDEs. Nevertheless, for bounded-domain problems involving the fractional Laplacian, three difficulties remain central: accurate evaluation of singular nonlocal operators, robust approximation of boundary-singular solutions, and scalability in high dimensions.
	
For high-dimensional fractional PDEs, Monte Carlo methods have become an important direction because they avoid the explicit construction of dense nonlocal discretization matrices. In the PINN framework, Guo et al. \cite{guo2022monte} introduced the MC-fPINNs method, in which the fractional Laplacian is approximated by Monte Carlo sampling after splitting the singular integral into a neighborhood \(B_{r_0}(\boldsymbol{x})\) of \(\boldsymbol{x}\) and its complement:
\begin{equation}\label{eq:FL_2part}
		(-\Delta)^{\alpha/2}u(\boldsymbol{x})
		=
		C_{d,\alpha}\left(
		\int_{\boldsymbol{y}\in B_{r_0}(\boldsymbol{x})}
		\frac{u(\boldsymbol{x})-u(\boldsymbol{y})}{\|\boldsymbol{x}-\boldsymbol{y}\|^{d+\alpha}}
		\,d\boldsymbol{y}
		+
		\int_{\boldsymbol{y}\notin B_{r_0}(\boldsymbol{x})}
		\frac{u(\boldsymbol{x})-u(\boldsymbol{y})}{\|\boldsymbol{x}-\boldsymbol{y}\|^{d+\alpha}}
		\,d\boldsymbol{y}
		\right).
\end{equation}
Here the first term is singular, whereas the second is regular. This splitting is dimension-friendly and naturally suited to the PINN framework. However, directly approximating the singular part by Monte Carlo sampling may introduce substantial variance and sampling error, which in turn limits the attainable accuracy and slows convergence, especially in very high dimensions. Moreover, the overall performance is often sensitive to user-specified parameters, such as the splitting radius and the truncation used in the singular part. These issues may become more pronounced for larger fractional orders, for which the kernel is more strongly singular and the numerical treatment is more delicate.
Several recent works have improved Monte Carlo methods for high-dimensional fractional PDEs from different perspectives. Hu et al. \cite{hu2024tackling} replaced the one-dimensional singular radial Monte Carlo integration by Gauss--Jacobi quadrature, thereby reducing variance and improving both accuracy and convergence in very high dimensions. Wang and Karniadakis \cite{wang2024gmc} developed a more general Monte Carlo PINN framework for fractional PDEs on irregular domains. From a different viewpoint, Sheng et al. \cite{sheng2023efficient} proposed an efficient probabilistic Monte Carlo solver based on Green-function/Feynman--Kac representations and walk-on-spheres-type ideas, which avoids direct discretization of the hypersingular operator.
Despite these advances, existing splitting-based MC-fPINN discretizations for bounded-domain fractional Laplacian problems \cite{guo2022monte,hu2024tackling} still rely largely on a fixed splitting radius \(r_0\) in \cref{eq:FL_2part}. For example, in \cite{hu2024tackling}, \(r_0\) is chosen as the diameter of the support of the problem. Although the influence of this parameter on accuracy has been examined \cite{guo2022monte}, its choice remains problem-dependent and is generally difficult to determine a priori. Near the boundary, such a fixed-radius splitting does not fully exploit the local geometry and may mix interior and exterior contributions in a non-adaptive manner. As a result, the attainable accuracy and robustness of these methods remain limited, particularly when boundary singularities and exterior Dirichlet constraints both play an important role.

Motivated by the above developments \cite{guo2022monte,hu2024tackling}, we combine boundary-adaptive region decomposition with deterministic radial quadrature and Monte Carlo angular sampling to obtain a more accurate and robust discretization for bounded-domain fractional Laplacian problems. Then we develop a quadrature-enhanced Monte Carlo fractional physics-informed neural network (QE-MC-fPINN) method. The main contributions of this work are as follows:
\begin{enumerate}
		\item We introduce a spatially varying radius \(r_0(\boldsymbol{x})\), defined by the distance from \(\boldsymbol{x}\) to \(\partial\Omega\), together with directional distance-to-boundary information, to decompose the fractional Laplacian into three parts: a near-field singular term, an interior far-field term, and an exterior far-field term. Gauss--Jacobi quadrature is used for the singular near-field radial integral, Gauss quadrature for the regular interior radial integral, and Monte Carlo sampling for the angular variables. This construction retains scalability in high dimensions while improving radial accuracy and boundary resolution.
		
		\item We embed the proposed discretization into a PINN trial space equipped with explicit feature functions, which improves training stability and enhances the approximation of low-regularity boundary-singular solutions under exterior Dirichlet constraints.
\end{enumerate}
	
The remainder of this paper is organized as follows. \cref{sec:Preli} introduces the model problems, the PINN formulation, the proposed Quadrature-Enhanced discretization, and related Monte Carlo discretizations. \cref{Sec:Numerical} presents the numerical results. Finally, \cref{sec:conclusion} concludes the paper.
	
\section{Problem Setting and PINN Formulation}\label{sec:Preli}
	We first focus on the fractional Poisson equation (fPE) with homogeneous Dirichlet exterior conditions.
	\begin{equation}\label{eq:fra_lapla_problem}
		\begin{aligned}
			(-\Delta)^{\alpha / 2} u(\boldsymbol{x}) &= f(\boldsymbol{x}), && \boldsymbol{x} \in \Omega, \\
			u(\boldsymbol{x}) &= 0, && \boldsymbol{x} \in \Omega^c
		\end{aligned}
	\end{equation}
defined on a bounded domain $\Omega \subset \mathbb{R}^d$.

We also consider the following fractional PDE on a bounded spatial domain $\Omega\subset\mathbb{R}^{d}$ \cite{hu2024tackling,pang2019fpinns}:
\begin{equation}\label{eq:fPDE}
		\begin{aligned}
			\frac{\partial^{\gamma}u(\boldsymbol{x},t)}{\partial t^\gamma}
			+ c(-\Delta)^{\alpha/2}u(\boldsymbol{x},t)
			+ \boldsymbol{v}\cdot\nabla_{\boldsymbol{x}}u(\boldsymbol{x},t)
			&= f(\boldsymbol{x},t),
			&& (\boldsymbol{x},t)\in\Omega\times(0,T], \\
			u(\boldsymbol{x},0)
			&= u_0(\boldsymbol{x}),
			&& \boldsymbol{x}\in\Omega, \\
			u(\boldsymbol{x},t)
			&= 0,
			&& (\boldsymbol{x},t)\in\Omega^{c}\times(0,T].
		\end{aligned}
\end{equation}
Here, $f(\boldsymbol{x},t)$ is the forcing term and $\boldsymbol{v}\in\mathbb{R}^{d}$ is the advection velocity.
The Caputo time-fractional derivative is defined by
	\begin{equation}\label{eq:time_fractional}
		\frac{\partial^{\gamma}f(t)}{\partial t^\gamma}
		:=
		\frac{1}{\Gamma(1-\gamma)}
		\int_{0}^{t}
		(t-\tau)^{-\gamma}
		\frac{\partial f(\tau)}{\partial \tau}\,d\tau,
		\qquad 0<\gamma<1.
	\end{equation}
	For convenience, we introduce the operator
	\begin{equation}\label{eq:operator_phi}
		\mathcal{L}^{\phi}[u(\boldsymbol{x},t)]
		:=
		\frac{\partial^{\gamma}u(\boldsymbol{x},t)}{\partial t^\gamma}
		+c (-\Delta)^{\alpha/2}u(\boldsymbol{x},t)
		+ \boldsymbol{v}\cdot\nabla_{\boldsymbol{x}}u(\boldsymbol{x},t),
	\end{equation}
	where $\phi=\{\gamma,\alpha,\boldsymbol{v},c\}$ denotes the collection of model parameters. For the fractional Poisson equation \cref{eq:fra_lapla_problem}, the operator reduces to \((-\Delta)^{\alpha/2}\).

\subsection{A Feature-enhanced PINN formulation}
	
In this work, we employ physics-informed neural networks (PINNs) \cite{raissi2019physics} to solve the forward problems associated with the time-dependent fractional PDE \cref{eq:fPDE} and the fractional Poisson equation \cref{eq:fra_lapla_problem}. These two problems can be handled in a similar manner within the same PINN framework. The basic idea is to approximate the unknown solution by a neural network and to determine the network parameters by minimizing a loss function constructed from the residual of the governing equation together with the associated initial and boundary conditions. As a result, the neural-network approximation is trained not only to fit collocation points, but also to satisfy the physical law prescribed by the PDE.
	
A well-known difficulty in PINN training is the imbalance between different components of the loss function, especially between the PDE residual term and the terms enforcing the initial and boundary conditions. This difficulty is more pronounced for fractional PDEs, because the nonlocal nature of the fractional Laplacian often leads to reduced boundary regularity and singular solution profiles near $\partial\Omega$. In particular, for fractional Dirichlet problems, the solution typically exhibits an $\frac{\alpha}{2}$-order boundary singularity of the form
	\begin{equation*}
		u(\boldsymbol{x}) \approx \operatorname{dist}(\boldsymbol{x},\partial\Omega)^{\alpha/2}u_{\mathrm{reg}}(\boldsymbol{x}),
		\qquad \boldsymbol{x}\in\overline{\Omega},
	\end{equation*}
where \(u_{\mathrm{reg}}\) is comparatively smoother \cite{Grubb2015FractionalLO,Guo2023ADL}. Therefore, capturing the boundary behavior accurately is essential in the design of effective neural-network trial spaces.
	
\begin{figure}[htbp]
		\centering
		\begin{tikzpicture}[
			node distance=5mm and 7mm,
			>=latex,
			every node/.style={font=\small},
			line/.style={->, thick},
			box/.style={
				draw,
				rounded corners,
				align=center,
				inner sep=2.5pt,
				minimum height=0.82cm
			},
			op/.style={
				draw,
				circle,
				inner sep=0.5pt,
				minimum size=5.3mm
			}
			]
			
			\node[box, minimum width=1.4cm] (input) {Input $(\boldsymbol{x},t)$};
			\node[box, minimum width=2.0cm, right=8mm of input] (phi) {NN basis\\
				$\varphi_j(\boldsymbol{x},t;\theta)$};
			\node[op, right=7mm of phi] (mul) {$\times$};
			\node[box, minimum width=2.8cm, right=7mm of mul] (psi) {
				$\psi_j=t^\gamma b(\boldsymbol{x})^{\mu_j}\varphi_j(\boldsymbol{x},t)$};
			\node[op, right=7mm of psi] (sum) {$\sum$};
			\node[box, minimum width=2.5cm, right=7mm of sum] (trial)
			{$\Psi=\sum_{j=1}^p \psi_j+u_0$};
			
			\node[box, minimum width=0.6cm, minimum height=0.4cm, above=4mm of mul] (time) {$t^\gamma$};
			\node[box, minimum width=0.6cm, minimum height=0.3cm, below=4mm of mul] (bdry) {$b(\boldsymbol{x})^{\mu_j}$};
			
			\node[box, minimum width=2.6cm, below=10mm of trial] (loss)
			{Residual loss \\ $\mathcal{L}_{\mathrm{res}}(\theta)$};
			\node[box, minimum width=1.8cm, left=7mm of loss] (opt)
			{$\min_{\theta}\mathcal{L}_{\mathrm{res}}$};
			
			\draw[line] (input) -- (phi);
			\draw[line] (phi) -- (mul);
			\draw[line] (time) -- (mul);
			\draw[line] (bdry) -- (mul);
			\draw[line] (mul) -- (psi);
			\draw[line] (psi) -- (sum);
			\draw[line] (sum) -- (trial);
			\draw[line] (trial.south) -- (loss.north);
			\draw[line] (loss.west) -- (opt.east);
			\draw[->, thick, dashed] (opt.west) -| ([yshift=-1mm]phi.south);
		\end{tikzpicture}
		\caption{Architecture of the proposed feature-enhanced PINN . The trial function incorporates the temporal factor \(t^\gamma\) and the boundary feature \(\left(b(\boldsymbol{x})\right)^{\mu_j}\) to enforce homogeneous constraints and improve the approximation of low-regularity solutions. The displayed trial function corresponds to the time-dependent case.}
		\label{fig:feature_enhanced_pinn}
	\end{figure}
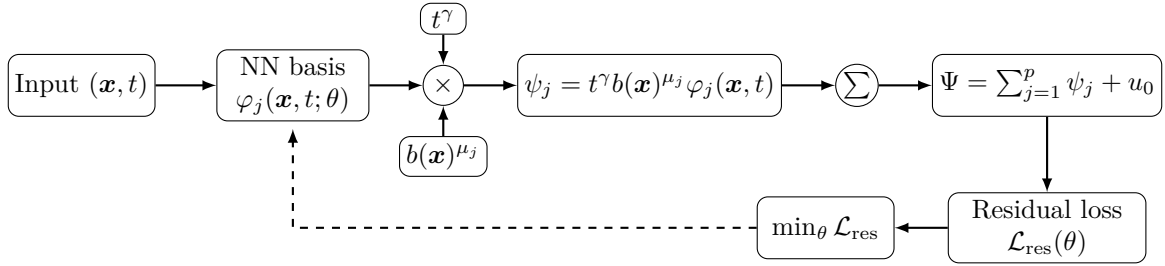
	
To impose the homogeneous boundary condition in a hard manner, we introduce a prescribed boundary feature function $b(\boldsymbol{x})$ satisfying
	\[
	b(\boldsymbol{x})>0,\qquad \boldsymbol{x}\in\Omega,
	\]
and
	\[
	b(\boldsymbol{x})=0,\qquad \boldsymbol{x}\in\Omega^{c}.
	\]
Let $\Psi(\boldsymbol{x},t;\theta)$ denote the neural-network approximation to the exact solution $u(\boldsymbol{x},t)$, where $\theta$ collects all trainable parameters. For the time-dependent problem \cref{eq:fPDE}, we define the trial solution by
	\begin{equation}\label{eq:pinn_trial_solution}
		\Psi(\boldsymbol{x},t;\theta)
		=
		u_0(\boldsymbol{x})
		+
		\sum_{j=1}^{p}\psi_j(\boldsymbol{x},t;\theta),
		\qquad
		\psi_j(\boldsymbol{x},t;\theta)
		=
		t^{\gamma}\left(b(\boldsymbol{x})\right)^{\mu_j}\varphi_j(\boldsymbol{x},t;\theta),
	\end{equation}
where $\{\varphi_j(\boldsymbol{x},t;\theta)\}_{j=1}^{p}$ are neural-network-generated basis functions.  For the fractional Poisson equation \cref{eq:fra_lapla_problem}, the corresponding trial space is obtained by omitting the temporal factor \(t^\gamma\) and the initial-condition term \(u_0(\boldsymbol{x})\).
The factor $t^\gamma$ is introduced to better capture the typical near-initial-time behavior associated with the Caputo fractional derivative  and to improve the approximation quality near $t=0$ \cite{lin2025TFPIDE}.
The exponents $\{\mu_j\}_{j=1}^{p}$ are boundary feature parameters associated with $b(\boldsymbol{x})$. They play an important role in representing possible singular behavior of the exact solution near the boundary, and hence significantly influence the approximation accuracy. 
We choose the exponents \(\{\mu_j\}_{j=1}^{p}\) by linear interpolation:
	\begin{equation}\label{eq:def_muj}
		\mu_j
		=
		\frac{\alpha}{2}+(j-1)\frac{\mu_p-\frac{\alpha}{2}}{p-1},
		\qquad	j=1,2,\dots,p.
	\end{equation}
Then \(\{\mu_j\}_{j=1}^{p}\) forms a monotonically increasing sequence from \(\frac{\alpha}{2}\) to \(\mu_p\). In particular, \(\mu_1=\frac{\alpha}{2}\) determines the leading boundary behavior of the trial function. This choice is inspired by the adaptive feature constructions in \cite{Guo2023ADL,lin2025TFPIDE}.
	
By construction, the trial function satisfies the homogeneous boundary condition exactly; for the time-dependent problem, it also enforces the initial condition exactly. Hence only the interior PDE residual needs to be minimized during training. Let
	\[
	\left\{(\boldsymbol{x}_{n}^{\mathrm{res}},t_{n}^{\mathrm{res}})\right\}_{n=1}^{N_{\mathrm{res}}}
	\subset \Omega\times(0,T]
	\]
	be the interior collocation points. The residual loss is defined by
	\begin{equation}\label{eq:loss_residual}
		\mathcal{L}_{\mathrm{res}}(\theta)
		=
		\frac{1}{N_{\mathrm{res}}}
		\sum_{n=1}^{N_{\mathrm{res}}}
		\left|
		\mathcal{L}^{\phi}\Psi\!\left(\boldsymbol{x}_{n}^{\mathrm{res}},t_{n}^{\mathrm{res}};\theta\right)
		-
		f\!\left(\boldsymbol{x}_{n}^{\mathrm{res}},t_{n}^{\mathrm{res}}\right)
		\right|^{2}.
	\end{equation}
The resulting training problem reads
	\begin{equation}\label{eq:pinn_minimization}
		\theta^{*}\in\arg\min_{\theta}\mathcal{L}_{\mathrm{res}}(\theta).
	\end{equation}
	
This yields a unified mesh-free framework for fractional PDEs. Its performance depends on the accurate evaluation of the residual, which in turn requires reliable discretizations of the nonlocal integral terms appearing in the fractional operator.
	
\subsection{Gauss--Jacobi Quadrature}\label{sec:G_J_qua}
	
Gauss--Jacobi quadrature is an efficient and high-precision method for approximating definite integrals. It is designed for integrals over $[-1,1]$ with the weight function $(1-x)^{\beta_{1}}(1+x)^{\beta_{2}}$ and takes the form (see \cite{Shen2011Spectral,brzezinski2018computation})
	\begin{equation}\label{Def_gauss_jacobi_11}
\int_{-1}^1 (1-x)^{\beta_1}(1+x)^{\beta_2} f(x)\,dx
\approx
\sum_{i=1}^N \widehat{w}_i^{(\beta_1,\beta_2)} f(\widehat{x}_i^{(\beta_1,\beta_2)}).
	\end{equation}
where $\beta_{1}> -1$ and $\beta_{2}> -1$.
	
The quadrature formula for a general interval $[a, b]$ is constructed as follows:
	\begin{equation}\label{Gauss_Jacobi_2}
		\int_a^b (b-x)^{\beta_{1}}(x-a)^{\beta_{2}} f(x)\,dx
		\approx
		\sum_{k=1}^{N} w_{k}^{(\beta_{1},\beta_{2})} f(x_{k}^{(\beta_{1},\beta_{2})}),
	\end{equation}
where the nodes and weights are obtained through the affine transformation
	\begin{equation*}
		x_{k}^{(\beta_{1},\beta_{2})} = \frac{b-a}{2} \widehat{x}_{k}^{(\beta_{1},\beta_{2})} + \frac{a+b}{2}, \quad
		w_{k}^{(\beta_{1},\beta_{2})} = \left( \frac{b-a}{2} \right)^{\beta_{1}+\beta_{2}+1} \widehat{w}_{k}^{(\beta_{1},\beta_{2})}.
	\end{equation*}
This quadrature rule will be repeatedly used for the singular radial integrals arising in the fractional Laplacian and for the weakly singular temporal integrals in the Caputo derivative. In particular, we will use the case \((\beta_1,\beta_2)=(0,1-\alpha)\) for singular radial integrals and \((-\gamma,\gamma-1)\) for the Caputo derivative after a suitable change of variables.
	
\subsection{Computation of the Fractional Laplacian}\label{sec:Computation_fractional_Laplacian}
In this subsection, we present the proposed discretization of the fractional Laplacian and briefly review related Monte Carlo discretizations.
\subsubsection{Quadrature-Enhanced MC-fPINN}
For notational simplicity, we write \(\psi_j(\boldsymbol{x},t)\) for \(\psi_j(\boldsymbol{x},t;\theta)\). We consider the fractional Laplacian defined in \cref{eq:fractional_laplacian}. For the trial function \(\Psi\) in \cref{eq:pinn_trial_solution}, one has
	\begin{equation*}
		(-\Delta)^{\alpha/2}\Psi(\boldsymbol{x},t;\theta)
		=
		\sum_{j=1}^{p}(-\Delta)^{\alpha/2}\psi_j(\boldsymbol{x},t)
		+
		(-\Delta)^{\alpha/2}u_0(\boldsymbol{x}).
	\end{equation*}

\begin{figure}[!htbp]
	\centering
	{\includegraphics[width=10cm]{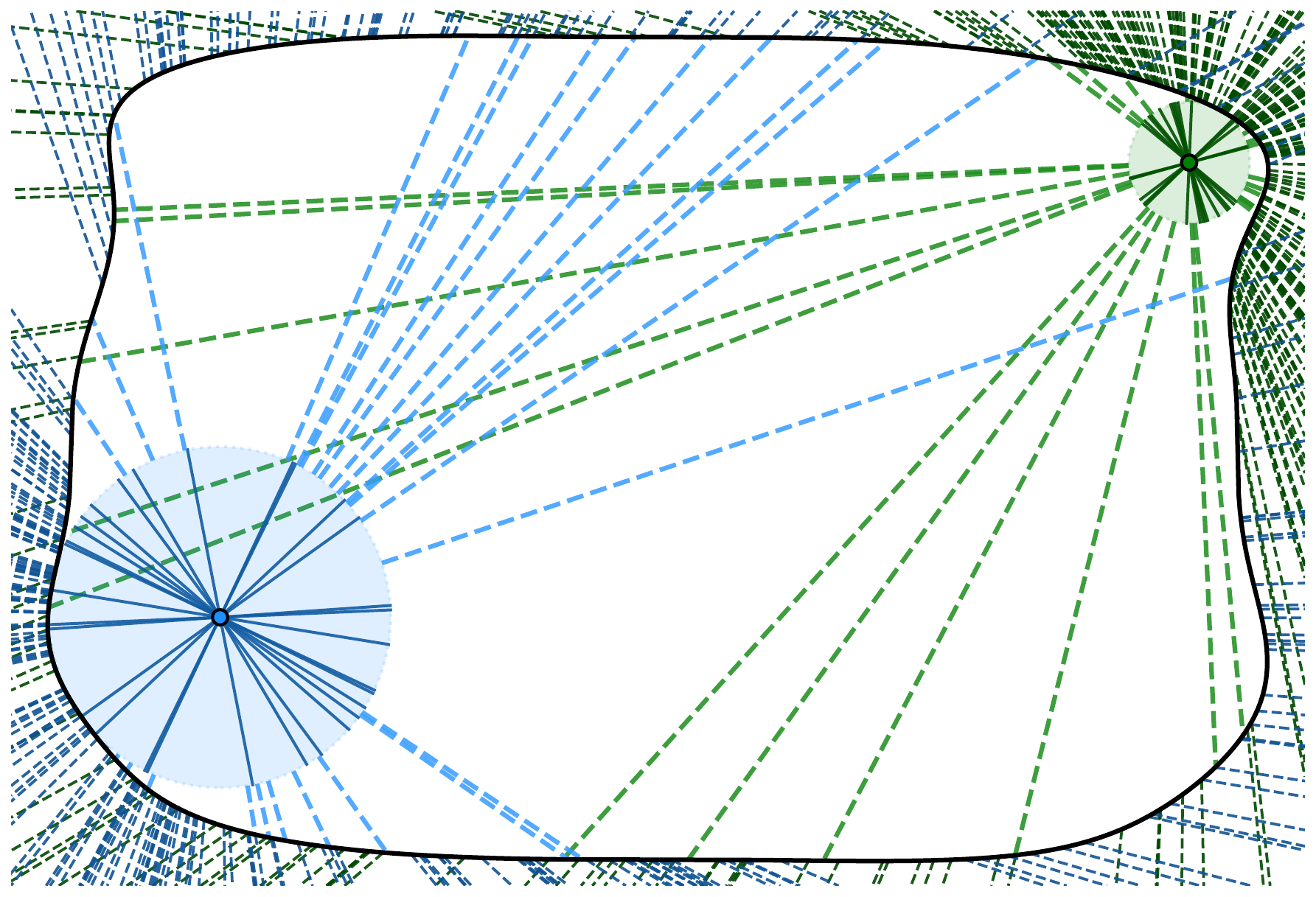}}
	\caption{Discretization scheme of the fractional Laplacian  on a 2D irregular domain. 
		The black closed curve denotes the physical boundary $\partial\Omega$. 
		For each interior evaluation point $\boldsymbol{x}$ (colored markers), a local ball $B_{r_0(\boldsymbol{x})}(\boldsymbol{x})$ of radius $r_0(\boldsymbol{x})$ is defined. 
		The plane $\mathbb{R}^2$ is partitioned into three regions and discretized using three independent sets of Monte Carlo angular directions:
		(1) $\boldsymbol{x} \to \partial B_{r_0(\boldsymbol{x})}(\boldsymbol{x})$ (near-field, 15 symmetrically distributed directions);
		(2) $\partial B_{r_0(\boldsymbol{x})}(\boldsymbol{x}) \to \partial\Omega$ (interior far-field, 30 directions);
		(3) $\partial\Omega \to \infty$ (exterior far-field, 256 directions).
		Ray segments represent directional quadrature paths, and colors distinguish different interior points.
	}
	\label{fig:2d_fra_dis}
\end{figure}
Our discretization uses deterministic quadrature in the radial direction and Monte Carlo sampling in the angular direction. Following the general philosophy of \cite{hu2024tackling}, this strategy is implemented here within a boundary-adaptive decomposition based on the spatially varying radius \(r_0(\boldsymbol{x})\), defined as the minimum distance from \(\boldsymbol{x}\) to \(\partial\Omega\), together with directional distance-to-boundary information. Accordingly, for each \(\psi_j\), we split the fractional Laplacian into near-field and far-field parts:
	\begin{equation*}
		\frac{(-\Delta)^{\alpha/2}\psi_j(\boldsymbol{x},t)}{C_{d,\alpha}}
		=
		\left(
		\int_{\|\boldsymbol{y}-\boldsymbol{x}\|<r_0(\boldsymbol{x})}
		+
		\int_{\|\boldsymbol{y}-\boldsymbol{x}\|\ge r_0(\boldsymbol{x})}
		\right)
		\frac{\psi_j(\boldsymbol{x},t)-\psi_j(\boldsymbol{y},t)}
		{\|\boldsymbol{x}-\boldsymbol{y}\|^{d+\alpha}}
		\,d\boldsymbol{y}
		=: I_{1,j}(\boldsymbol{x},t)+I_{2,j}(\boldsymbol{x},t).
	\end{equation*}
This yields a three-part decomposition of the operator into singular near-field, interior far-field, and exterior far-field contributions. A schematic illustration in a two-dimensional irregular domain is shown in \cref{fig:2d_fra_dis}.  For each interior point \(\boldsymbol{x}\), the local radius \(r_0(\boldsymbol{x})\) defines the near-field region, while the remaining domain is split into interior and exterior far-field parts.

The near-field integral is singular at $\boldsymbol{y}=\boldsymbol{x}$. To handle this singular behavior accurately, we employ the Gauss--Jacobi quadrature introduced in \cref{sec:G_J_qua}, which provides an efficient and high-precision numerical treatment \cite{jahanshahi2017fractional,lin2025TFPIDE,wang2024gmc}.
The near-field contribution is given by
	\begin{equation*}
		I_{1,j}(\boldsymbol{x},t)
		=
		r_0(\boldsymbol{x})^{-\alpha}
		\int_{\mathcal{S}_{+}^{d-1}}
		\int_0^1
		\frac{
			F_j\bigl(\boldsymbol{x},t,r_0(\boldsymbol{x})\tau,\boldsymbol{\xi}\bigr)
		}{
			\tau^{1+\alpha}
		}
		\,d\tau\,d\sigma(\boldsymbol{\xi}),
	\end{equation*}
where
	\begin{equation*}
		F_j(\boldsymbol{x},t,r,\boldsymbol{\xi})
		:=
		2\psi_j(\boldsymbol{x},t)
		-\psi_j(\boldsymbol{x}+r\boldsymbol{\xi},t)
		-\psi_j(\boldsymbol{x}-r\boldsymbol{\xi},t),
	\end{equation*}
and $d\sigma(\boldsymbol{\xi})$ denotes the surface measure on the unit sphere. Rewriting the radial integral in weighted form, we obtain
	\begin{equation*}
		I_{1,j}(\boldsymbol{x},t)
		=
		r_0(\boldsymbol{x})^{-\alpha}
		\int_{\mathcal{S}_{+}^{d-1}}
		\int_0^1
		\tau^{1-\alpha}
		\frac{
			F_j\bigl(\boldsymbol{x},t,r_0(\boldsymbol{x})\tau,\boldsymbol{\xi}\bigr)
		}{
			\tau^2
		}
		\,d\tau\,d\sigma(\boldsymbol{\xi}).
	\end{equation*}
We then approximate the radial part by Gauss--Jacobi quadrature and the angular part by Monte Carlo sampling:
	\begin{equation}\label{eq:near_singular}
		I_{1,j}(\boldsymbol{x},t)
		\approx
		\frac{|\mathcal{S}^{d-1}|\,r_0(\boldsymbol{x})^{-\alpha}}{M'}
		\sum_{\ell=1}^{M'}
		\sum_{k=1}^{N}
		w_k^{(0,1-\alpha)}
		\frac{
			F_j\!\left(
			\boldsymbol{x},t,
			r_0(\boldsymbol{x})\tau_k^{(0,1-\alpha)},
			\boldsymbol{\xi}_{\ell}
			\right)
		}{
			\bigl(\tau_k^{(0,1-\alpha)}\bigr)^2
		}.
	\end{equation}
Here, \(\{\boldsymbol{\xi}_{\ell}\}_{\ell=1}^{M'}\) are uniformly distributed random directions on the upper hemisphere \(\mathcal{S}_{+}^{d-1}\). By symmetry, the integral over the upper hemisphere is equivalent to one half of the full-sphere integral, which leads to the prefactor \(|\mathcal S^{d-1}|\) in \cref{eq:near_singular}.

\begin{remark}
	For the near-field contribution, the accuracy of the discretization is determined by the behavior of \(F_j(\boldsymbol{x},t,r,\boldsymbol{\xi})\) as \(r\to0\).
	Assuming \(\psi_j\) is smooth enough, we expand \(\psi_j(\boldsymbol{x} \pm r\boldsymbol{\xi},t)\) in Taylor series:
	\begin{equation*}
		\psi_j(\boldsymbol{x} \pm r\boldsymbol{\xi},t) = \psi_j(\boldsymbol{x},t)\pm r \nabla \psi_j(\boldsymbol{x},t) \cdot \boldsymbol{\xi} + \frac{r^{2}}{2} \boldsymbol{\xi}^{\top} H(\psi_j)(\boldsymbol{x},t) \boldsymbol{\xi}\pm \frac{r^3}{6} D^3\psi_j(\boldsymbol{x},t)(\boldsymbol{\xi},\boldsymbol{\xi},\boldsymbol{\xi}) + O(r^4),
	\end{equation*}
	where \(H(\psi_j)(\boldsymbol{x},t)\) is the Hessian matrix of \(\psi_j\) at \(\boldsymbol{x}\), 
	containing all second-order partial derivatives and \(D^3\psi_j(\boldsymbol{x},t)(\boldsymbol{\xi},\boldsymbol{\xi},\boldsymbol{\xi})\) 
	denotes the third-order directional derivative of \(\psi_j\) at \(\boldsymbol{x}\) 
	in the direction \(\boldsymbol{\xi}\).
	Substituting these into the definition of \(F_j/r^{1+\alpha}\) yields
	\begin{equation}\label{eq:taylor}
		\frac{F_j(\boldsymbol{x},t, r, \boldsymbol{\xi})}{r^{1+\alpha}} = - r^{1-\alpha} \boldsymbol{\xi}^{\top} H(\psi_j)(\boldsymbol{x},t) \boldsymbol{\xi} + O(r^{3-\alpha}), \quad r \to 0^+,
	\end{equation}
	which reveals that the integrand exhibits the asymptotic 
	behavior \(r^{1-\alpha}\) near \(r=0\). For \(\alpha\in (0,2)\), this singularity is 
	integrable, i.e., weak singularity. After integration, we have 
	\[
	\int_0^{r_{0}} \frac{F_j(\boldsymbol{x},t, r, \boldsymbol{\xi})}{r^{1+\alpha}} \, dr = -\frac{r_{0}^{2-\alpha}}{2-\alpha} \boldsymbol{\xi}^{\top} H(\psi_j)(\boldsymbol{x},t) \boldsymbol{\xi} + O(r_{0}^{4-\alpha}).
	\]
	This asymptotic structure motivates the use of Gauss--Jacobi quadrature 
	with parameters \((0, 1-\alpha)\), which is specifically 
	designed to handle  singularities of the form \(r^{1-\alpha}\). For any point \(\boldsymbol{x}\in\Omega\), since \(B_{r_0(\boldsymbol{x})}(\boldsymbol{x})\subset\Omega\), the local Taylor expansion is valid for \(0<r\le r_0(\boldsymbol{x})\) provided that \(\psi_j\) is sufficiently smooth in a neighborhood of \(\boldsymbol{x}\).  This justifies the use of Gauss--Jacobi quadrature and explains its high accuracy for approximating the radial integral over \((0,r_0(\boldsymbol{x})]\).
\end{remark}
	
The far-field contribution is decomposed as
	\begin{equation*}
		I_{2,j}(\boldsymbol{x},t)
		=
		\int_{\mathcal{S}^{d-1}}
		\bigl[
		Q_{1,j}(\boldsymbol{x},t,\boldsymbol{\xi})
		+
		Q_{2,j}(\boldsymbol{x},t,\boldsymbol{\xi})
		\bigr]
		\,d\sigma(\boldsymbol{\xi}),
	\end{equation*}
where
	\begin{equation*}
		Q_{1,j}(\boldsymbol{x},t,\boldsymbol{\xi})
		=
		\int_{r_0(\boldsymbol{x})}^{d_{\boldsymbol{x}}(\boldsymbol{\xi})}
		\frac{
			\psi_j(\boldsymbol{x},t)
			-
			\psi_j(\boldsymbol{x}+r\boldsymbol{\xi},t)
		}{
			r^{1+\alpha}
		}
		\,dr
	\end{equation*}
denotes the interior far-field contribution and
	\begin{equation*}
		d_{\boldsymbol{x}}(\boldsymbol{\xi})
		:=
		\min\{s>0:\boldsymbol{x}+s\boldsymbol{\xi}\in\partial\Omega\}
	\end{equation*}
is the distance from $\boldsymbol{x}$ to the boundary along the direction $\boldsymbol{\xi}$. By using the variable transformation \(r=r(\boldsymbol{x},\boldsymbol{\xi },s)=r_{0}(\boldsymbol{x})+\bigl( d_{\boldsymbol{x}}(\boldsymbol{\xi })-r_{0}(\boldsymbol{x}) \bigr)s\) we have
	\begin{align*}
		Q_{1,j}(\boldsymbol{x},t,\boldsymbol{\xi })
		&=\int_{r_{0}(\boldsymbol{x})}^{d_{\boldsymbol{x}}(\boldsymbol{\xi })}
		\frac{\psi_j(\boldsymbol{x},t) - \psi_j(\boldsymbol{x}+r\boldsymbol{\xi},t)}{r^{1+\alpha}}dr\\
		&=\int_0^{1}{\frac{\psi_j(\boldsymbol{x},t) - \psi_j(\boldsymbol{x}
				+r(\boldsymbol{x},\boldsymbol{\xi },s)\boldsymbol{\xi},t)}{r(\boldsymbol{x},\boldsymbol{\xi },s)^{1+\alpha}}\left(d_{\boldsymbol{x}}(\boldsymbol{\xi})-r_{0}(\boldsymbol{x})\right)}ds\\
		&\approx \left(d_{\boldsymbol{x}}(\boldsymbol{\xi})-r_{0}(\boldsymbol{x})\right)
		\sum_{m=1}^{N_{0}}{w_{m}\frac{\psi_j(\boldsymbol{x},t) - \psi_j(\boldsymbol{x}
				+r(\boldsymbol{x},
				\boldsymbol{\xi },s_{m})\boldsymbol{\xi},t)}{r(\boldsymbol{x},\boldsymbol{\xi },s_{m})^{1+\alpha}}},
	\end{align*}
where \(N_0\) denotes the number of Gauss quadrature points on \([0,1]\), and \(\{s_m\}_{m=1}^{N_0}\), \(\{w_m\}_{m=1}^{N_0}\) are the corresponding nodes and weights.
The exterior contribution is given analytically by
	\begin{equation*}
		Q_{2,j}(\boldsymbol{x},t,\boldsymbol{\xi})
		=
		\psi_j(\boldsymbol{x},t)
		\int_{d_{\boldsymbol{x}}(\boldsymbol{\xi})}^{\infty}
		\frac{dr}{r^{1+\alpha}}
		=
		\frac{1}{\alpha}\,d_{\boldsymbol{x}}(\boldsymbol{\xi})^{-\alpha}\psi_j(\boldsymbol{x},t).
	\end{equation*}


The radial integral in $Q_{1,j}$ is still evaluated by Gauss quadrature, while the angular integral is approximated by Monte Carlo sampling. Consequently, the far-field term is approximated by
	\begin{equation*}
		I_{2,j}(\boldsymbol{x},t)
		\approx
		\sum_{i=1}^{2}
		\sum_{k=1}^{M_i}
		\frac{|\mathcal{S}^{d-1}|}{M_i}
		\,Q_{i,j}(\boldsymbol{x},t,\boldsymbol{\xi}_{ik}),
	\end{equation*}
where $\{\boldsymbol{\xi}_{ik}\}_{k=1}^{M_i}$ are uniformly distributed random directions on $\mathcal{S}^{d-1}$. Here, $M_1$ denotes the number of Monte Carlo samples used for the interior far-field term $Q_{1,j}$, and $M_2$ denotes the number of Monte Carlo samples used for the exterior term $Q_{2,j}$.

Combining the above approximations, we obtain a hybrid Monte Carlo discretization of the fractional Laplacian. 
In this scheme, the singular near-field radial integral is approximated by Gauss--Jacobi quadrature, the interior far-field radial integral is approximated by Gauss quadrature, and the angular integrations are approximated by Monte Carlo sampling. This construction preserves the high accuracy of deterministic quadrature in the radial direction while retaining the flexibility and scalability of Monte Carlo sampling in the angular direction. Moreover, owing to the spatially varying radius $r_0(\boldsymbol{x})$ and the decomposition of the far-field region into interior and exterior parts, the present Monte Carlo discretization is more adaptive to the local boundary geometry than the improved MC-fPINN method in \cite{hu2024tackling}. Therefore, it yields a more accurate and robust discretization for bounded-domain problems involving the fractional Laplacian.
\begin{remark}
	The present discretization differs from the method in \cite{hu2024tackling} in three main respects.
	First, in the near-field discretization \cref{eq:near_singular}, we exploit the symmetry of the integrand and restrict Monte Carlo sampling to the upper hemisphere. This reduces the number of angular evaluations by one half without loss of accuracy, and hence improves computational efficiency.
	Second, the use of the spatially varying radius \(r_0(\boldsymbol{x})\) leads to a three-part decomposition of the fractional Laplacian, in which the regular interior far-field contribution appears as a separate term. This decomposition is more adaptive to the local geometry of the domain and yields a more accurate discretization on bounded domains.
	Third, the analytical exterior contribution \(Q_{2,j}\) retains directional information, which allows the associated spherical integral to be approximated with a denser set of directions when higher accuracy is needed. This provides greater flexibility for high-accuracy angular approximation of the exterior contribution than in the fixed-radius formulation of \cite{hu2024tackling}. Consequently, the present formulation can yield a more accurate discretization of the fractional Laplacian on bounded domains, especially when geometric effects are significant.
\end{remark}
\subsubsection{Related Monte Carlo discretizations}
\paragraph{MC-fPINN}
	
Guo et al.~\cite{guo2022monte} introduced the Monte Carlo fractional physics-informed neural network (MC-fPINN) for discretizing nonlocal operators. For the fractional Laplacian, the regional decomposition yields a singular near-field part and a regular far-field part, each of which can be approximated via Monte Carlo sampling. The corresponding discretization can be written as
	\begin{align*}
		(-\Delta)^{\alpha/2}u(\boldsymbol{x})&=	C_{d,\alpha}\,\frac{|\mathcal{S}^{d-1}|\,r_0^{\,2-\alpha}}{2\left(2-\alpha\right)}\,
		\mathbb{E}_{\boldsymbol{\xi},\,r_I\sim f_I(r_I)}\left[
		\frac{2u(\boldsymbol{x})-u(\boldsymbol{x}-r_\epsilon\,\boldsymbol{\xi})
			-u(\boldsymbol{x}+r_\epsilon\,\boldsymbol{\xi})}{r_\epsilon^2}
		\right]\\
		&+	C_{d,\alpha}\,\frac{|\mathcal{S}^{d-1}|\,r_0^{-\alpha}}{2\alpha}\,
		\mathbb{E}_{\boldsymbol{\xi},\,r_o\sim f_O(r_o)}\left[
		2u(\boldsymbol{x})-u(\boldsymbol{x}-r_o\,\boldsymbol{\xi})
		-u(\boldsymbol{x}+r_o\,\boldsymbol{\xi})
		\right],
	\end{align*}
	where $\boldsymbol{\xi}$ is uniformly distributed on the unit $(d-1)$-sphere $\mathcal{S}^{d-1}$. 
	For the near-field contribution, the radial random variable $r_I$ follows the density
	\[
	f_I(r_I)
	=
	\frac{2-\alpha}{r_0^{\,2-\alpha}}\;r_I^{\,1-\alpha}\,
	\mathbf{1}_{r\in[0,r_0]},
	\]
	which can be sampled by the transformation
	$r_I/r_0\sim\operatorname{Beta}(2-\alpha,1)$. And $r_\epsilon=\max \{\epsilon,r_I \}$, where $\epsilon > 0$ is a small positive number.
	Similarly, for the far-field part, the radial variable follows
	\[
	f_O(r_o)
	=
	\alpha\,r_0^{\,\alpha}\,r_o^{-1-\alpha}\,
	\mathbf{1}_{r_o\in[r_0,\infty)},
	\]
	which can be sampled via
	$r_0/r_o\sim\operatorname{Beta}(\alpha,1)$.
	
	\paragraph{Improved MC-fPINN}
	
	Hu et al.~\cite{hu2024tackling} considered a comparable regional decomposition for the fractional Laplacian. The decomposition yields
	\[
	(-\Delta)^{\alpha/2}u(\boldsymbol{x})
	=\frac{ C_{d,\alpha}}{2}\int_{\mathcal{S}^{d-1}}
	\int_0^{r_0}\frac{
		2u(\boldsymbol{x})-u(\boldsymbol{x}-r\,\boldsymbol{\xi})
		-u(\boldsymbol{x}+r\,\boldsymbol{\xi})
	}{r^{1+\alpha}}\,dr\,d\boldsymbol{\xi}
	+
	  C_{d,\alpha}\frac{\lvert \mathcal{S}^{d-1}\rvert\, r_0^{-\alpha}}{\alpha}\,u(\boldsymbol{x}).
	\]

	To discretize the near-field part, they replace Monte Carlo sampling of the radial variable by a Gauss--Jacobi quadrature rule:
	\[
	\frac{|\mathcal{S}^{d-1}|}{2}
	\sum_{k=1}^{N}
	w_k^{(0,1-\alpha)}
	\frac{
		2u(\boldsymbol{x})-u(\boldsymbol{x}-\tau_k^{(0,1-\alpha)}\,\boldsymbol{\xi}_k)
		-u(\boldsymbol{x}+\tau_k^{(0,1-\alpha)}\,\boldsymbol{\xi}_k)
	}{
		\bigl(\tau_k^{(0,1-\alpha)}\bigr)^2
	}.
	\]
Compared with the MC-fPINN method, this strategy improves the radial discretization of the singular part and leads to better accuracy and faster convergence in high dimensions \cite{hu2024tackling}.

\section{Numerical examples}\label{Sec:Numerical}
In this section, we provide several examples to validate the efficiency and accuracy of 
the proposed neural network-based machine learning method.  
All the experiments are conducted on NVIDIA GeForce RTX 4090 D GPUs.

The following relative $\ell^2$ error is used to measure the accuracy of the approximate solution $\Psi(\boldsymbol{x},t;\theta)$ with respect to the exact solution $u(\boldsymbol{x},t)$:
\begin{equation*}
	e_{\mathrm{test}}
	:=
	\frac{
		\left(
		\sum\limits_{k=1}^{K}
		\bigl(
		\Psi(\boldsymbol{x}^{k},t^{k};\theta)-u(\boldsymbol{x}^{k},t^{k})
		\bigr)^2
		\right)^{1/2}
	}{
		\left(
		\sum\limits_{k=1}^{K}
		\bigl(
		u(\boldsymbol{x}^{k},t^{k})
		\bigr)^2
		\right)^{1/2}
	}.
\end{equation*}
Here, \(\{(\boldsymbol{x}^{k},t^{k})\}_{k=1}^{K}\subset \Omega\times(0,T]\) denotes the set of test points. 
In all computations, we use \(N_\tau=8\) Gauss--Jacobi quadrature points for the Caputo time-fractional derivative and set \(\mu_p=\mu_1+1\). 

For memory efficiency in high dimensions, the Monte Carlo sums are evaluated in chunks; this affects only the order of accumulation and does not change the underlying discretization. 

For brevity, in the tables and figures only, MC\(^{*}\), I-MC\(^{*}\), and QE-MC\(^{*}\) are used to denote MC-fPINN, Improved MC-fPINN, and Quadrature-Enhanced MC-fPINN, respectively.
\begin{table}[!htbp]
	\centering
	\caption{Exact solutions in the unit ball and their fractional Laplacians.}
	\label{tab:unit_ball_exact}
	\begin{tabular}{cc}
			\toprule
			$u(\boldsymbol{x})$ in the unit ball & $(-\Delta)^{\alpha/2}u(\boldsymbol{x})$ in the unit ball  \cite{guo2022monte,hu2024tackling,pang2019fpinns} \\
			\midrule
			$\left(1-\|\boldsymbol{x}\|^2\right)^{\alpha/2}$
			&
			$-2^\alpha \Gamma\!\left(\frac{\alpha}{2}+1\right)
			\Gamma\!\left(\frac{\alpha+d}{2}\right)
			\Gamma\!\left(\frac{d}{2}\right)^{-1}$
			\\[1.2ex]
			
			$\left(1-\|\boldsymbol{x}\|^2\right)^{1+\alpha/2}$
			&
			$-2^\alpha \Gamma\!\left(\frac{\alpha}{2}+2\right)
			\Gamma\!\left(\frac{\alpha+d}{2}\right)
			\Gamma\!\left(\frac{d}{2}\right)^{-1}
			\left(1-\left(1+\frac{\alpha}{d}\right)\|\boldsymbol{x}\|^2\right)$
			\\[1.2ex]
			
			$\left(1-\|\boldsymbol{x}\|^2\right)^{\alpha/2}x_d$
			&
			$-2^\alpha \Gamma\!\left(\frac{\alpha}{2}+1\right)
			\Gamma\!\left(\frac{\alpha+d}{2}+1\right)
			\Gamma\!\left(\frac{d}{2}+1\right)^{-1}x_d$
			\\[1.2ex]
			
			$\left(1-\|\boldsymbol{x}\|^2\right)^{1+\alpha/2}x_d$
			&
			$-2^\alpha \Gamma\!\left(\frac{\alpha}{2}+2\right)
			\Gamma\!\left(\frac{\alpha+d}{2}+1\right)
			\Gamma\!\left(\frac{d}{2}+1\right)^{-1}
			\left(1-\left(1+\frac{\alpha}{d+2}\right)\|\boldsymbol{x}\|^2\right)x_d$
			\\
			\bottomrule
		\end{tabular}
\end{table}

\subsection{Fractional Poisson equation}

Consider the fractional Poisson equation
\begin{equation}\label{eq:fPE}
	(-\Delta)^{\alpha/2}u(\boldsymbol{x})=f(\boldsymbol{x})
	\qquad \text{in } B_1=\{\boldsymbol{x}\in\mathbb{R}^d:\|\boldsymbol{x}\|<1\}.
\end{equation}
Several exact solutions and the corresponding source terms are listed in \cref{tab:unit_ball_exact}. To further examine anisotropic effects, we also consider the exact solution
\begin{equation}\label{eq:com_sol}
	u(\boldsymbol{x})
	=
	\left(1-\|\boldsymbol{x}\|^{2}\right)^{\alpha/2}
	\left(c_{1,0}+\sum_{i=1}^d c_{1,i}x_i\right)
	+
	\left(1-\|\boldsymbol{x}\|^{2}\right)^{1+\alpha/2}
	\left(c_{2,0}+\sum_{i=1}^d c_{2,i}x_i\right),
\end{equation}
where the coefficients \(c_{1,i}\) and \(c_{2,i}\) are sampled independently from the standard normal distribution \(\mathcal{N}(0,1)\) \cite{hu2024tackling}. The first term in \cref{eq:com_sol} contains the factor \((1-\|\boldsymbol{x}\|^2)^{\alpha/2}\), which reflects the leading boundary behavior typically observed in fractional Dirichlet problems on bounded domains. Consequently, this example is substantially more challenging than the smoother benchmark solutions in \cref{tab:unit_ball_exact}. Moreover, the affine factors introduce anisotropy, so that \cref{eq:com_sol} provides a stringent test of both boundary resolution and high-dimensional anisotropic approximation.

For a fair comparison among MC-fPINN, Improved MC-fPINN, and Quadrature-Enhanced MC-fPINN, we use the same PINN backbone and training protocol for all three methods. Specifically, all methods employ the same trial-solution ansatz with boundary factor
\[
b(\boldsymbol{x})=\max\bigl(1-\|\boldsymbol{x}\|^2,\,0\bigr),
\]
the same multilayer perceptron architecture with width \(128\) and depth \(4\), the same Tanh activation, and the same basis size \(p=16\). The numbers of residual and test points are fixed as \(100\) and \(20{,}000\), respectively, and all computations are carried out in double precision. The network parameters are optimized by Adam for \(20{,}000\) epochs with a learning-rate decay schedule. In addition, the number of Monte Carlo angular samples used for the near-field integral is fixed at \(64\) for all three methods. For the method-specific discretizations, MC-fPINN uses \(r_0=0.25\) and \(\varepsilon=10^{-6}\); Improved MC-fPINN uses \(r_0=2\), with the remaining settings taken from \cite{hu2024tackling}; and Quadrature-Enhanced MC-fPINN uses \(8\) Gauss--Jacobi points and \(10\) Gauss points for the radial discretization of the fractional Laplacian, together with \(64\) and \(256\) Monte Carlo angular samples for the interior far-field and analytical exterior far-field terms, respectively.  For the unit ball \(B_1\), the directional distance from an interior point \(\boldsymbol{x}\) to the boundary \(\partial B_1\) along a unit vector \(\boldsymbol{\xi}\) is
\[
d_{\boldsymbol{x}}(\boldsymbol{\xi})
=
-\boldsymbol{x}\cdot\boldsymbol{\xi}
+\sqrt{(\boldsymbol{x}\cdot\boldsymbol{\xi})^2+1-\|\boldsymbol{x}\|^2},
\qquad \|\boldsymbol{\xi}\|=1.
\]

\begin{table}[!htb]
	\centering
	\caption{Relative $\ell^2$ errors $e_{\mathrm{test}}$ for solving fractional Poisson equation (fPE) by MC-fPINN, Improved MC-fPINN \cite{hu2024tackling}, and Quadrature-Enhanced MC-fPINN.}
	\label{tab:ball_error}
	\small
	\begin{tabular}{cccccccc}
		\toprule
		\multirow{2}{*}{$d$} & \multirow{2}{*}{$\alpha$}
		& \multicolumn{3}{c}{$u(\boldsymbol{x})=(1-\|\boldsymbol{x}\|^{2})^{1+\alpha/2}$}
		& \multicolumn{3}{c}{$u(\boldsymbol{x})$ defined in \cref{eq:com_sol}} \\
		\cmidrule(lr){3-5} \cmidrule(lr){6-8}
		& & $\mathrm{MC}^{*}$ & I-$\mathrm{MC}^{*}$ & QE-$\mathrm{MC}^{*}$ & $\mathrm{MC}^{*}$ & I-$\mathrm{MC}^{*}$ & QE-$\mathrm{MC}^{*}$ \\
		\midrule
		\multirow{3}{*}{3}
		& 0.2 & 6.04e-3 & 2.04e-3 & 2.70e-4 & 1.23e-1  & 1.02e-1  & 1.07e-3 \\
		& 1.5  & 4.37e-2 & 2.89e-2 & 1.51e-3  & 4.10e-1  & 4.41e-1  & 1.49e-2 \\
		& 1.9  & 7.55e-1  & 6.42e-1 & 7.49e-4 & 4.63e-1 & 7.45e-1  & 6.15e-3 \\
		\cdashline{1-8}[1pt/1pt]
		\multirow{3}{*}{100}
		& 0.2 & 1.16e-2  & 9.00e-3 & 3.52e-4 & 2.75e-1  & 3.18e-1  & 2.74e-3 \\
		& 1.5  & 6.11e-2  & 1.00e-1  & 1.44e-4 & 6.02e-2 & 8.84e-2  & 1.35e-3 \\
		& 1.9  & 1.67e-1  & 2.67e-1  & 1.22e-4 & 2.34e-1  & 6.82e-1  & 2.46e-3 \\
		\cdashline{1-8}[1pt/1pt]
		\multirow{3}{*}{300}
		& 0.2 & 1.26e-2  & 1.11e-2   & 3.29e-4  & 5.01e-1 & 2.14e-1  & 2.37e-3 \\
		& 1.5  & 6.85e-2 & 9.78e-2 & 1.16e-4   & 2.38e-1 & 5.72e-2 & 4.89e-3 \\
		& 1.9  & 9.69e-2 & 4.20e-1 & 8.83e-5  & 4.63e-1 & 4.04e-1 & 5.43e-3 \\
		\cdashline{1-8}[1pt/1pt]
		\multirow{3}{*}{1,000}
		& 0.2 & 1.42e-2  & 1.28e-2  & 3.99e-4 & 2.56e-1  & 3.29e-1 & 1.88e-3 \\
		& 1.5  & 7.39e-2 & 9.54e-2 & 1.44e-4  & 8.94e-2 & 1.33e-1  & 3.96e-3 \\
		& 1.9  & 9.42e-2 & 2.31e-1  & 1.23e-4   & 1.07e-1  & 5.79e-1  & 3.76e-3 \\
		\cdashline{1-8}[1pt/1pt]
		\multirow{3}{*}{10,000}
		& 0.2 & 1.55e-2 & 1.50e-2  & 2.72e-4   & 3.93e-1  &  2.27e-1 & 4.24e-3\\
		& 1.5  & 9.17e-2 & 9.87e-2  & 9.17e-5   &  1.45e-1 & 3.78e-2  & 3.90e-3 \\
		& 1.9  & 1.12e-1  & 4.84e-1  & 7.00e-5   & 8.73e-2 & 2.19e-1  &  4.11e-3 \\
		\bottomrule
	\end{tabular}
\end{table}

\begin{figure}[!htbp]
	\centering
	\subfloat[$u(\boldsymbol{x})=(1-\|\boldsymbol{x}\|^{2})^{1+\alpha/2}$ , $d=1,000, \alpha =1.5$]{\includegraphics[width=0.48\textwidth]{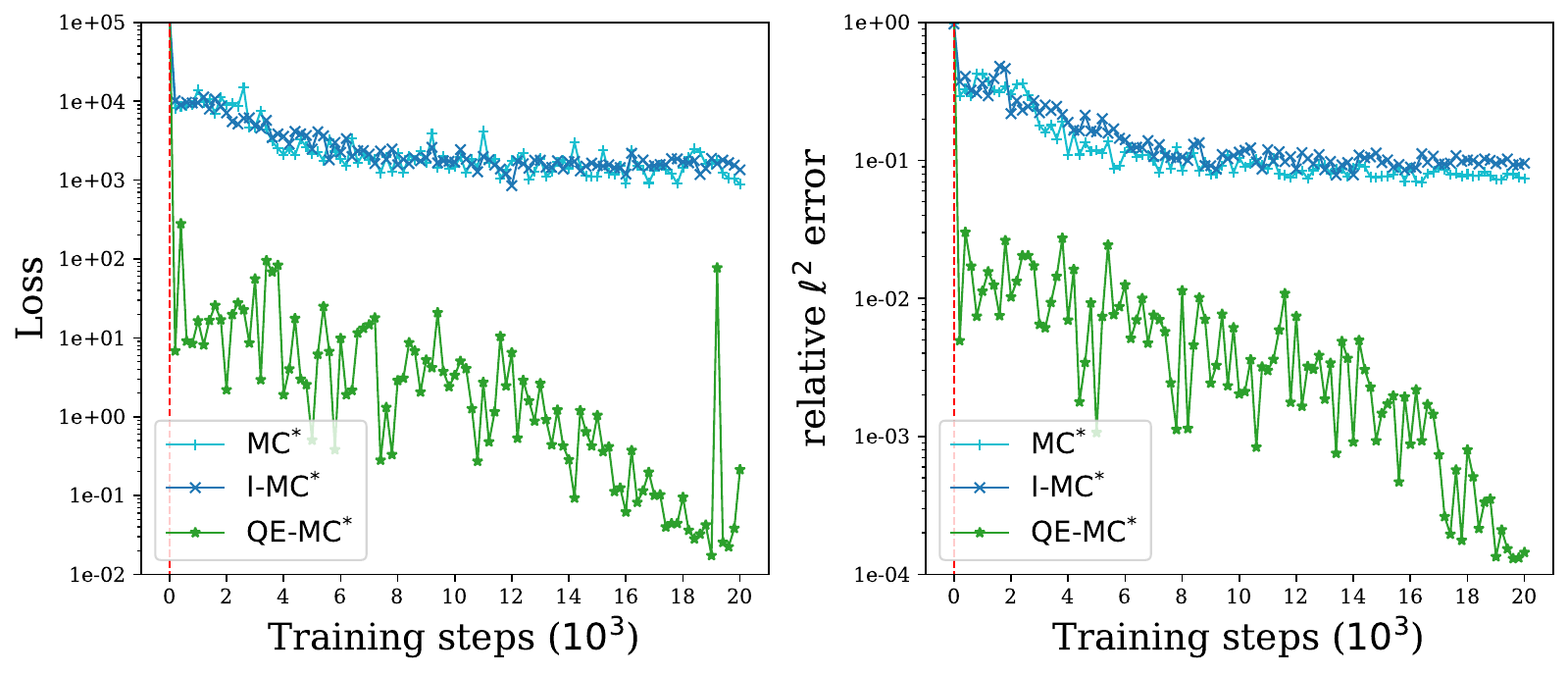}}
	\hfill
	\subfloat[$u(\boldsymbol{x})$ in \cref{eq:com_sol} , $d=1,000, \alpha =1.9$]{\includegraphics[width=0.48\textwidth]{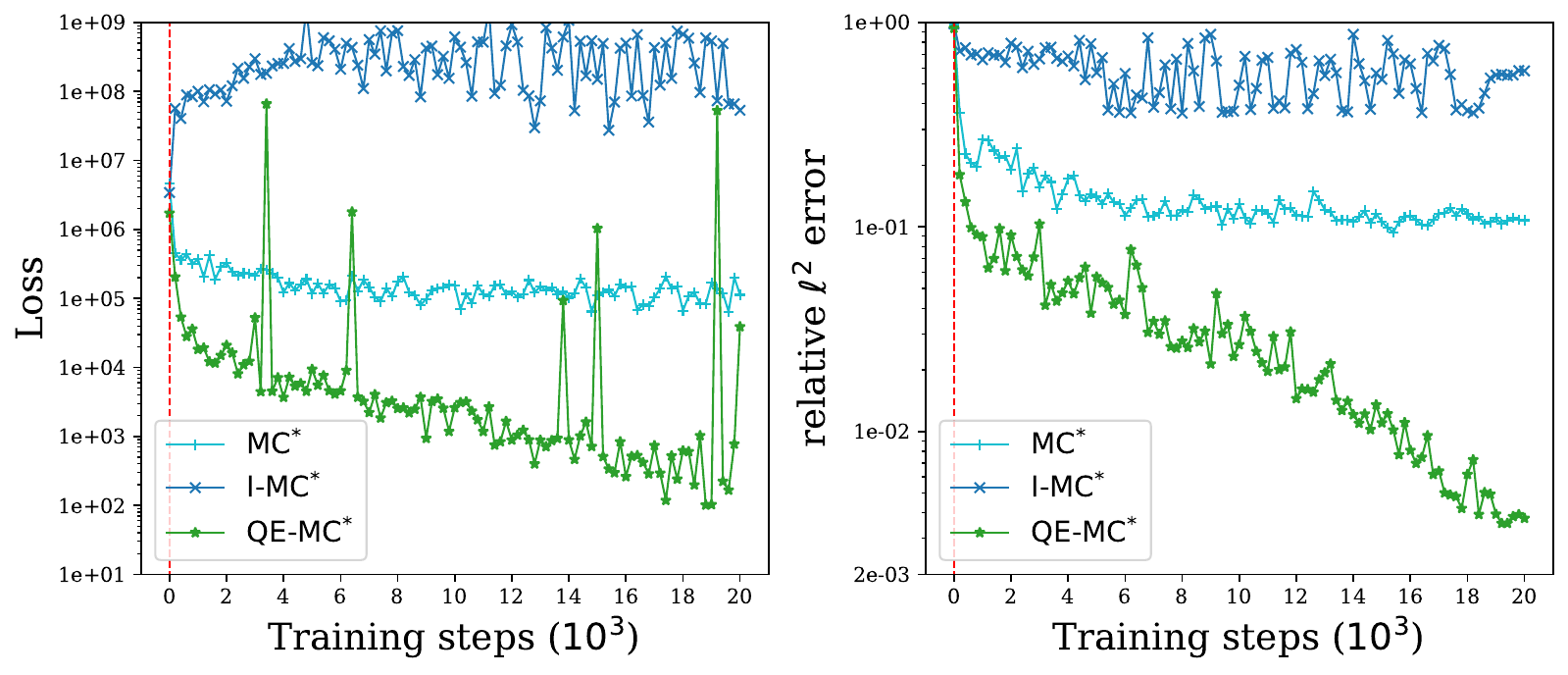}}
	
	\medskip  
	\subfloat[$u(\boldsymbol{x})=(1-\|\boldsymbol{x}\|^{2})^{1+\alpha/2}$ , $d=10,000, \alpha =1.5$]{\includegraphics[width=0.48\textwidth]{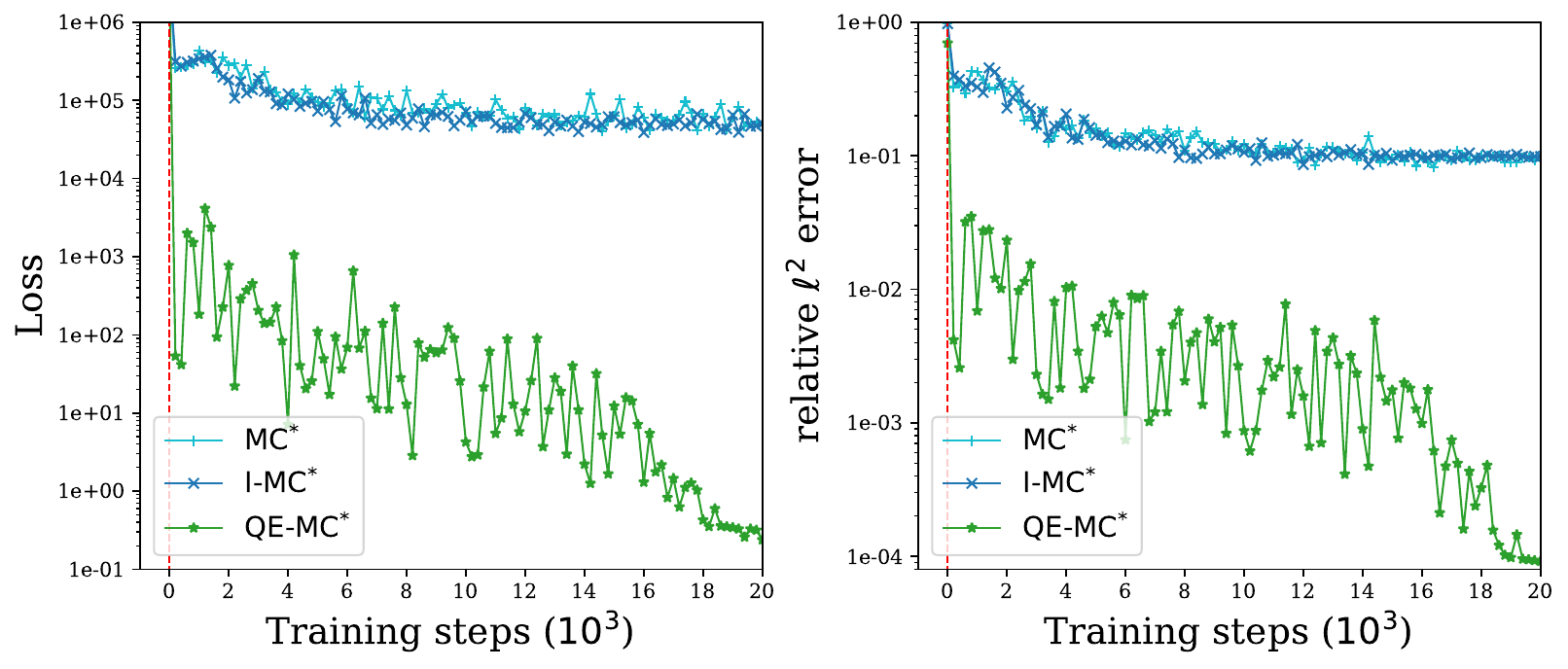}}
	\hfill
	\subfloat[$u(\boldsymbol{x})$  in \cref{eq:com_sol} , $d=10,000, \alpha =1.9$]{\includegraphics[width=0.48\textwidth]{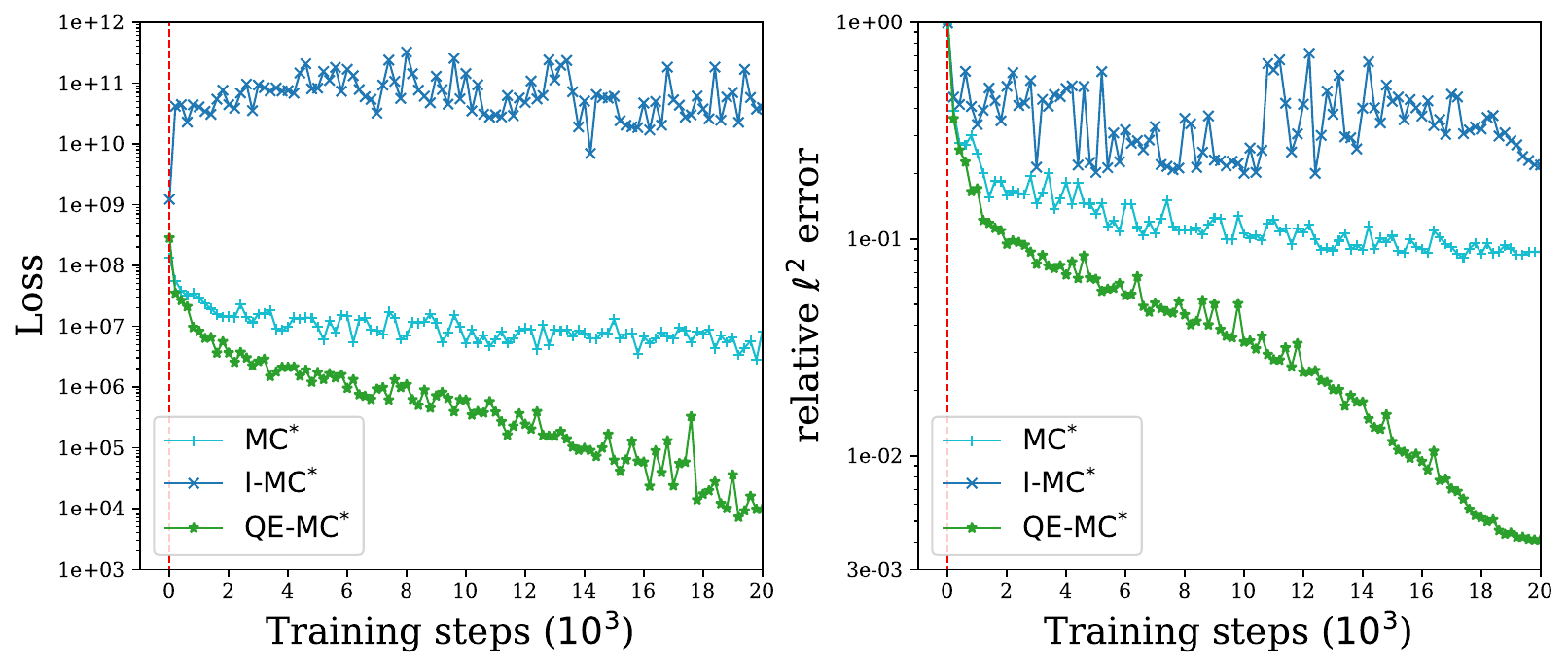}}
	
	\caption{Loss and relative $\ell^2$ error curves for the high-dimensional  fractional Poisson equation \cref{eq:fPE} on the unit ball.}
	\label{fig:loss_error_fPE}
\end{figure}

\cref{tab:ball_error,fig:loss_error_fPE} demonstrate that the proposed Quadrature-Enhanced MC-fPINN consistently delivers the best accuracy among the three methods over all tested dimensions and fractional orders. For the smooth benchmark solution \(u=(1-\|\boldsymbol{x}\|^2)^{1+\alpha/2}\), the proposed method reduces the relative error by one to several orders of magnitude compared with both MC-fPINN and Improved MC-fPINN, and this advantage becomes more pronounced as the dimension increases. For the more challenging anisotropic solution defined in \cref{eq:com_sol}, the improvement is even more significant, especially for larger fractional orders such as \(\alpha=1.5\) and \(\alpha=1.9\), where the two baseline methods often exhibit substantial loss of accuracy. The convergence histories in \cref{fig:loss_error_fPE} further show that Quadrature-Enhanced MC-fPINN attains a lower loss level and a smaller relative error with a more stable decay behavior in both \(d=1{,}000\) and \(d=10{,}000\). These results show that the proposed Quadrature-Enhanced decomposition substantially improves the robustness and resolution of the fractional Laplacian discretization in high dimensions.

\subsection{Time-Dependent Fractional Diffusion Equation}

The trial solution \(\Psi\) is defined in \cref{eq:pinn_trial_solution}. Since the initial condition \(u_0(\boldsymbol{x})\) is independent of time, its Caputo derivative vanishes, and thus
\[
\frac{\partial^{\gamma}}{\partial t^\gamma}\Psi(\boldsymbol{x},t;\theta)
=
\sum_{j=1}^{p}
b(\boldsymbol{x})^{\mu_j}
\,\partial_t^\gamma
\!\left(t^\gamma \varphi_j(\boldsymbol{x},t;\theta)\right).
\]

To discretize the Caputo time-fractional derivative, we employ the Gauss--Jacobi quadrature rule; see \cite{lin2025TFPIDE} for related details. For each basis term, we have
\begin{align}
	\frac{\partial^{\gamma}}{\partial t^\gamma}
	\Bigl(t^\gamma \varphi_j(\boldsymbol{x},t;\theta)\Bigr)
	&=
	\frac{1}{\Gamma(1-\gamma)}
	\int_0^t
	(t-s)^{-\gamma}
	\frac{\partial}{\partial s}
	\Bigl(s^\gamma \varphi_j(\boldsymbol{x},s;\theta)\Bigr)\,ds \nonumber\\
	&=
	\frac{1}{\Gamma(1-\gamma)}
	\int_0^1
	(1-\tau)^{-\gamma}\tau^{\gamma-1}
	S_{1,j}(\boldsymbol{x},t,\tau)\,d\tau,
	\qquad s=t\tau, \nonumber\\
	&\approx
	\frac{1}{\Gamma(1-\gamma)}
	\sum_{k=1}^{N_\tau}
	w_k^{(-\gamma,\gamma-1)}
	S_{1,j}\!\left(\boldsymbol{x},t,\tau_k^{(-\gamma,\gamma-1)}\right),
	\label{eq:GJ_t_gamma_basis}
\end{align}
where
$
S_{1,j}(\boldsymbol{x},t,\tau)
=
\gamma\,\varphi_j(\boldsymbol{x},t\tau;\theta)
+
t\tau\,\partial_t\varphi_j(\boldsymbol{x},t\tau;\theta).
$

We consider two exact solutions with different boundary regularity:
\begin{align}
	\label{eq:time_sol_smooth}
	u(\boldsymbol{x},t)
	&=
	t^{2.5}(1-\|\boldsymbol{x}\|^2)^{1+\alpha/2}
	\left(c_{2,0}+\sum_{i=1}^d c_{2,i}x_i\right),\\
	\label{eq:time_sol_singular}
	u(\boldsymbol{x},t)
	&=
	t^{2.5}(1-\|\boldsymbol{x}\|^2)^{\alpha/2}
	\left[
	\left(c_{1,0}+\sum_{i=1}^d c_{1,i}x_i\right)
	+
	(1-\|\boldsymbol{x}\|^2)
	\left(c_{2,0}+\sum_{i=1}^d c_{2,i}x_i\right)
	\right].
\end{align}
The first solution, \cref{eq:time_sol_smooth}, contains the higher-order boundary factor \((1-\|\boldsymbol{x}\|^2)^{1+\alpha/2}\) and is therefore relatively smooth near \(\partial B_1\). The second solution, \cref{eq:time_sol_singular}, contains the lower-order factor \((1-\|\boldsymbol{x}\|^2)^{\alpha/2}\), which reflects the leading boundary behavior typically associated with fractional Dirichlet problems on bounded domains. Consequently, \cref{eq:time_sol_singular} is more singular near the boundary and is numerically more challenging. In both examples, the affine factors introduce anisotropy, while the solutions vanish on the boundary of the unit ball. The coefficients \(c_{1,i}\) and \(c_{2,i}\) are sampled independently from the standard normal distribution \(\mathcal{N}(0,1)\).

For the numerical experiments for \cref{eq:fPDE}, we set
\[
\phi=\{\gamma,\alpha,\boldsymbol{v},c\}=\{0.5,1.5,1,1\},\qquad T=1.
\]
As indicated by the results in \cref{tab:ball_error}, MC-fPINN and Improved MC-fPINN already provide reasonable accuracy for the spatial discretization of this problem.

To ensure a fair comparison, we use the same PINN backbone and training protocol as in the fractional Poisson experiments. In particular, the spatial discretization of the fractional Laplacian is exactly the same as that used for \cref{eq:fPE}. The only additional discretization in the time-dependent case is that of the Caputo time-fractional derivative, for which we use \(N_\tau=8\) Gauss--Jacobi quadrature points. Unless otherwise stated, the numbers of residual and test points are fixed as \(100\) and \(20{,}000\), respectively; all computations are performed in double-precision floating-point arithmetic; and the network parameters are optimized by Adam for \(100{,}000\) epochs with a learning-rate decay schedule.
\begin{figure}[!h]
	\centering
	\subfloat[$u(\boldsymbol{x},t)$  in \cref{eq:time_sol_smooth}  with $d=3$]{\includegraphics[width=0.48\textwidth]{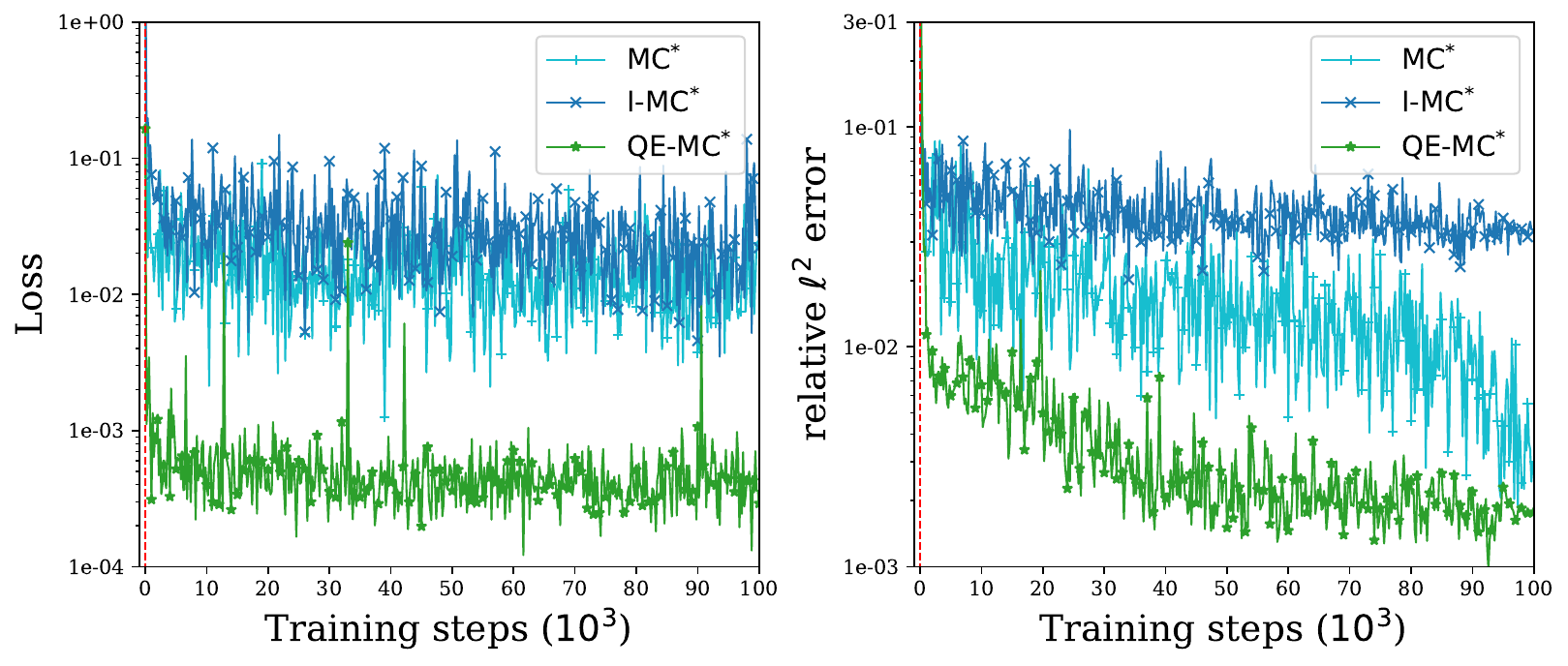}}
	\hfill
	\subfloat[$u(\boldsymbol{x},t)$  in \cref{eq:time_sol_singular} with $d=3$]{\includegraphics[width=0.48\textwidth]{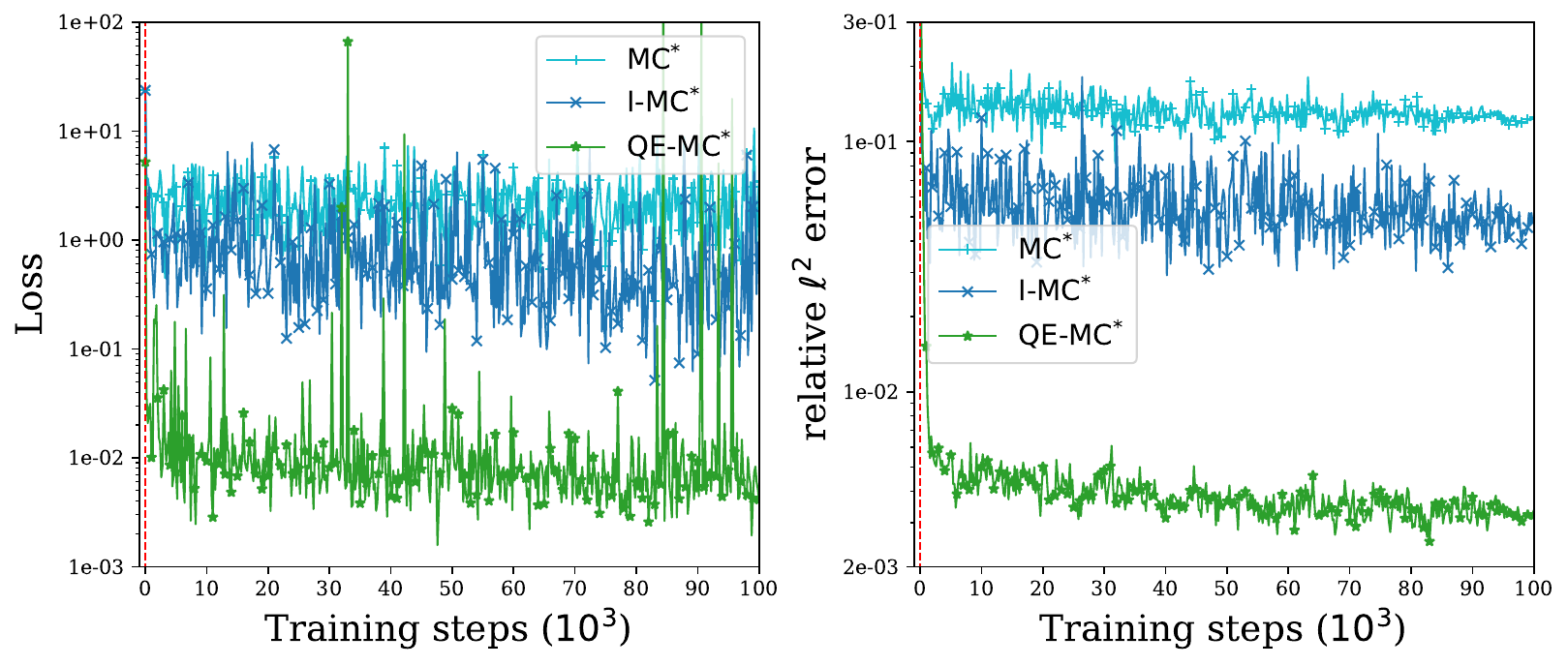}}
	
	\medskip  
	\subfloat[$u(\boldsymbol{x},t)$  in \cref{eq:time_sol_smooth}  with $d=100$]{\includegraphics[width=0.48\textwidth]{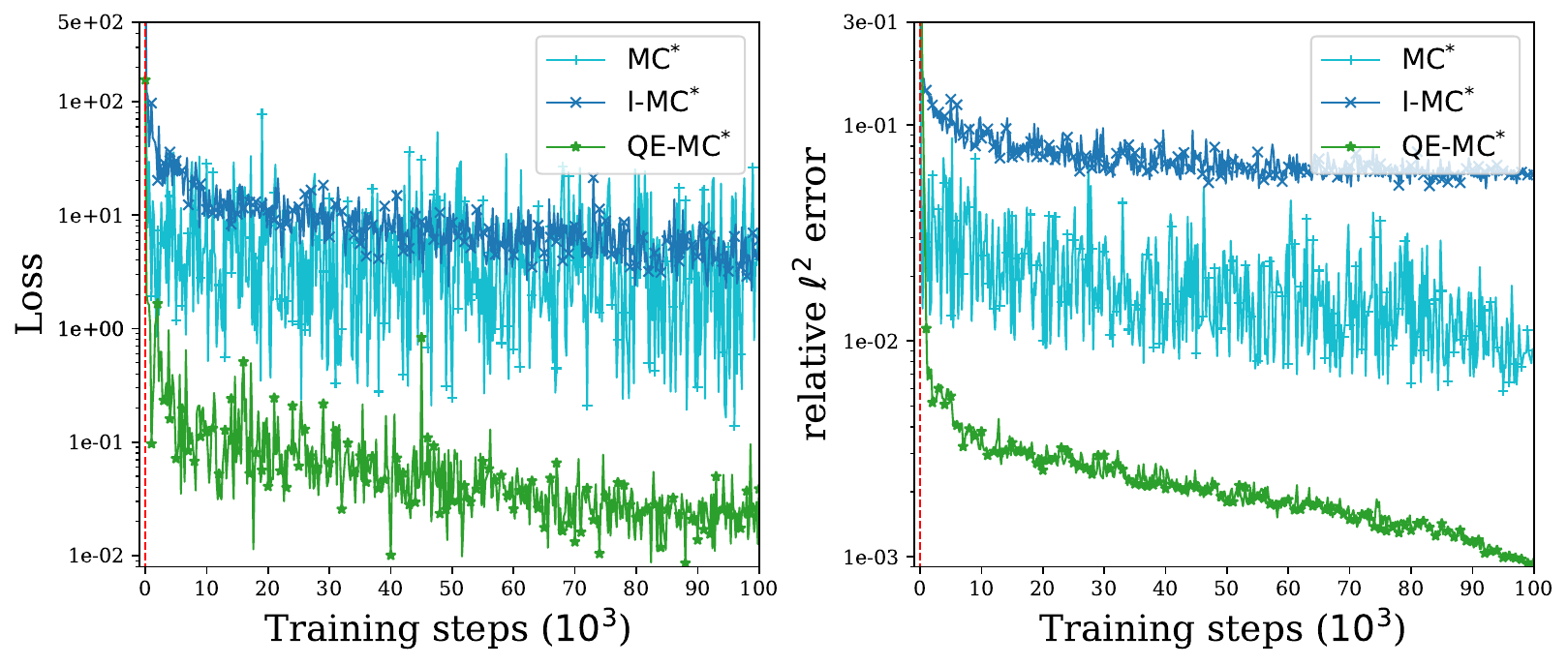}}
	\hfill
	\subfloat[$u(\boldsymbol{x},t)$  in \cref{eq:time_sol_singular} with $d=100$]{\includegraphics[width=0.48\textwidth]{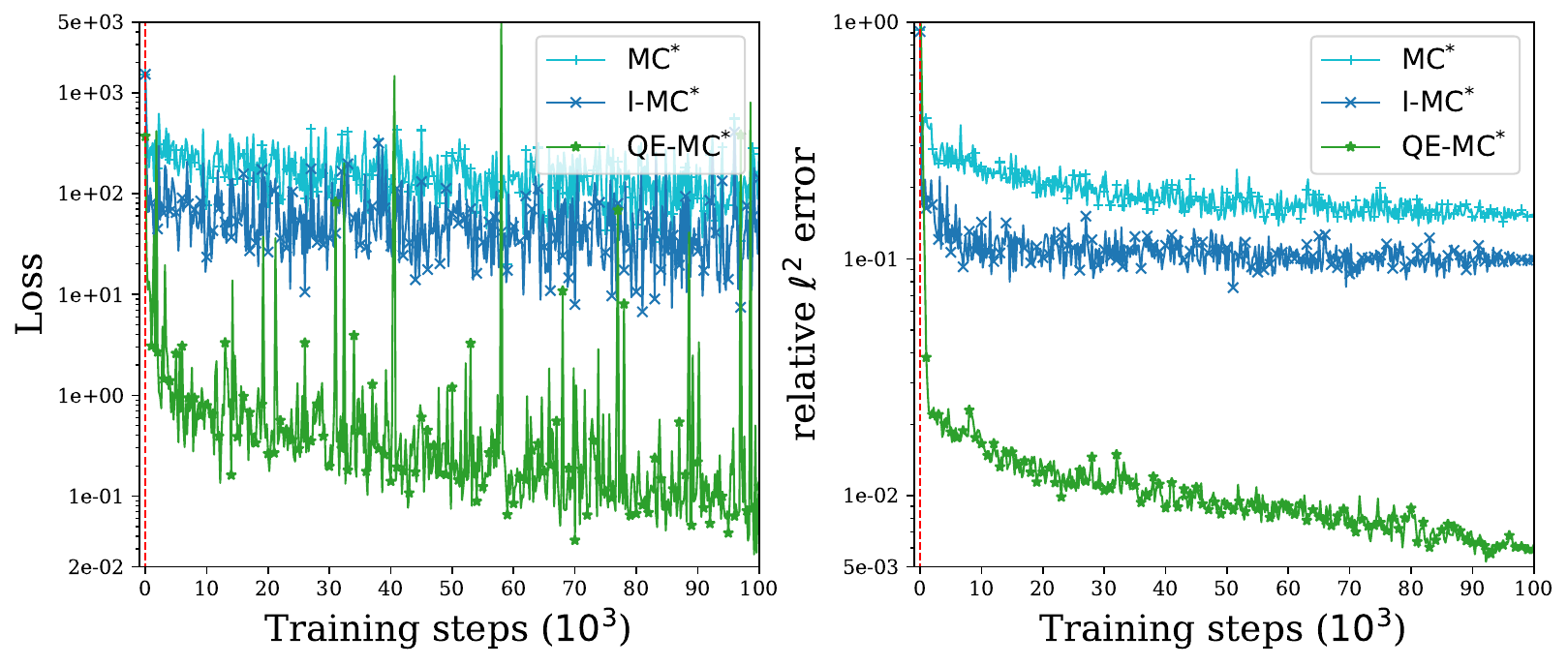}}
	
		\medskip  
	\subfloat[$u(\boldsymbol{x},t)$  in \cref{eq:time_sol_smooth}  with $d=1000$]{\includegraphics[width=0.48\textwidth]{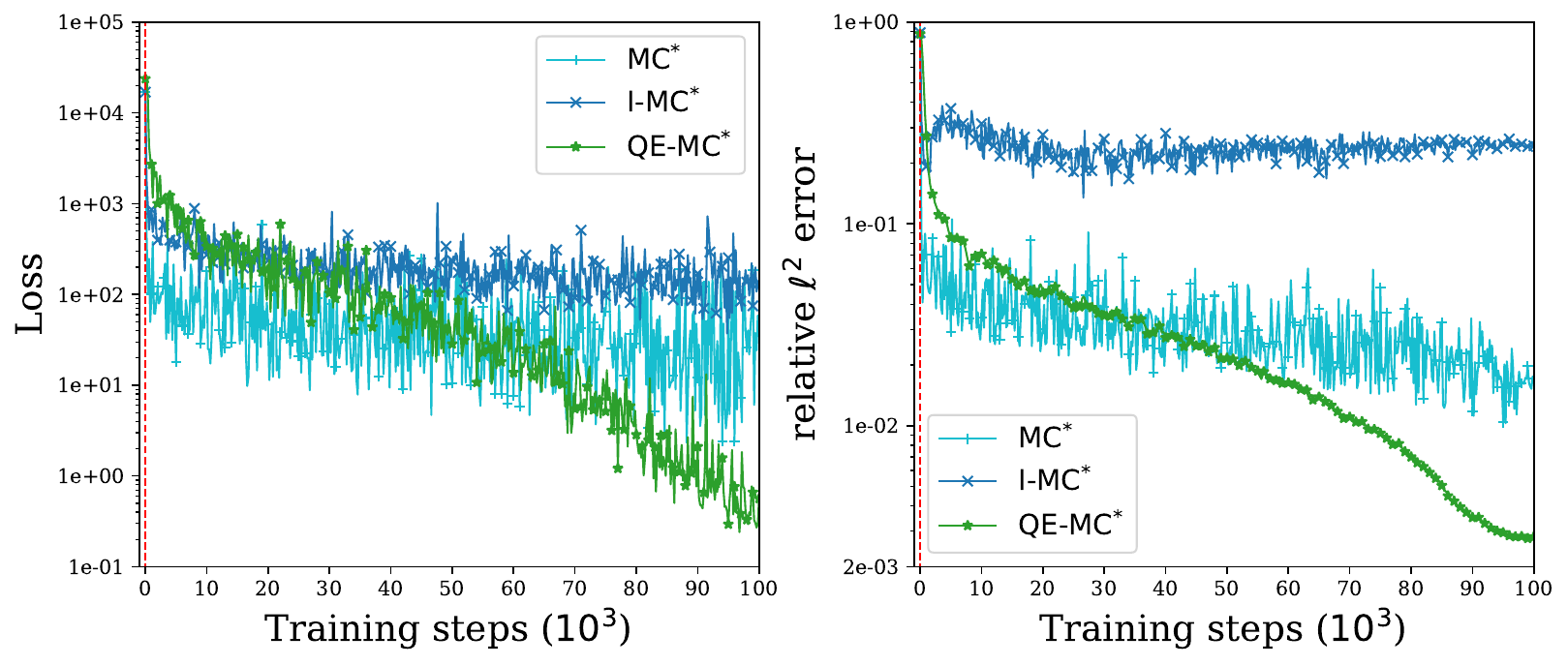}}
	\hfill
	\subfloat[$u(\boldsymbol{x},t)$  in \cref{eq:time_sol_singular} with $d=1000$]{\includegraphics[width=0.48\textwidth]{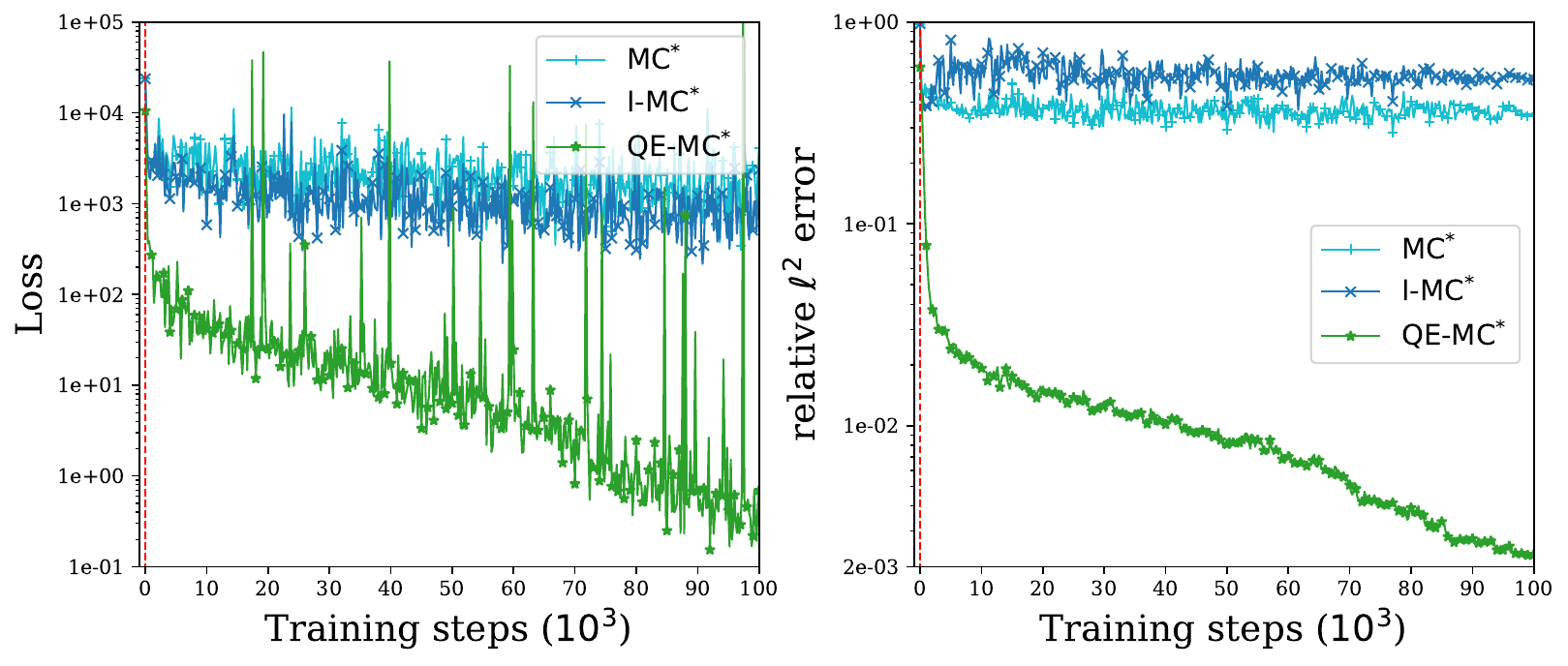}}
		
	\medskip  
	\subfloat[$u(\boldsymbol{x},t)$  in \cref{eq:time_sol_smooth}  with $d=5000$]{\includegraphics[width=0.48\textwidth]{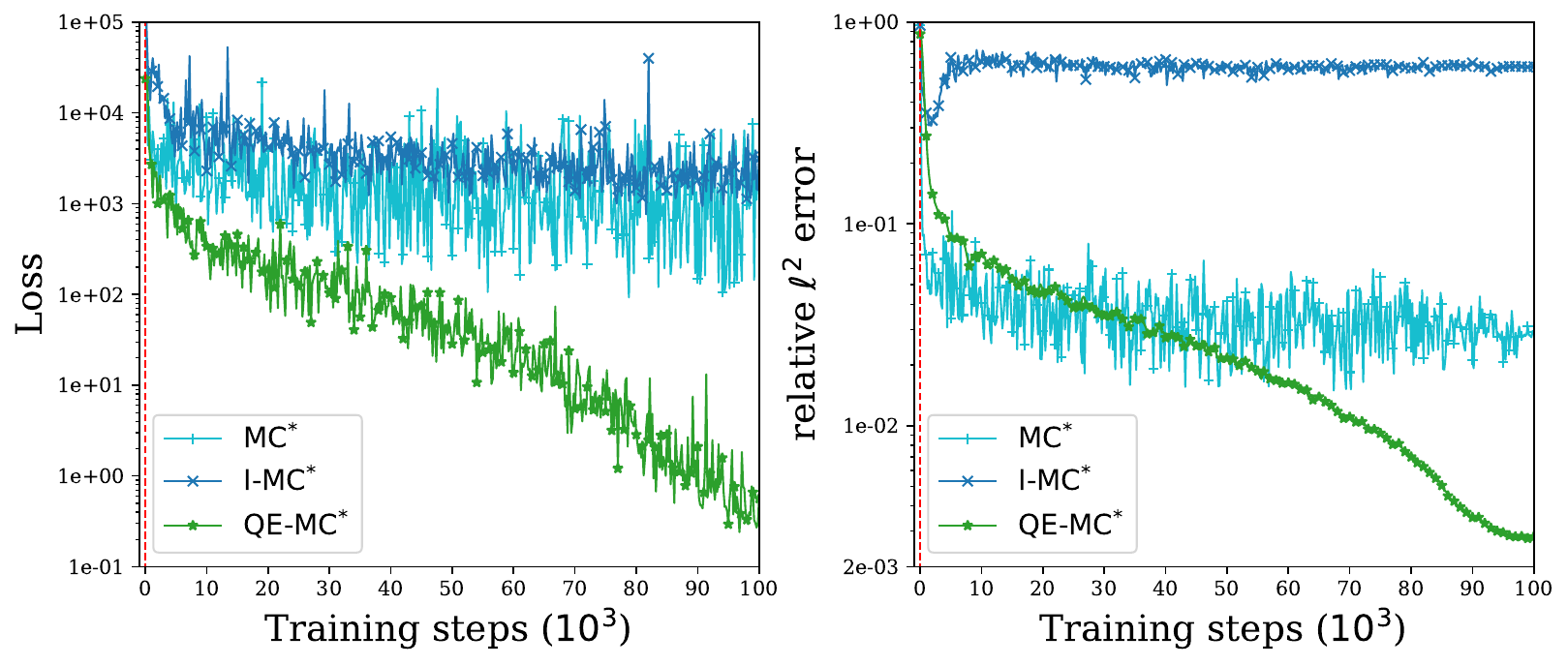}}
	\hfill
	\subfloat[$u(\boldsymbol{x},t)$  in \cref{eq:time_sol_singular} with $d=5000$]{\includegraphics[width=0.48\textwidth]{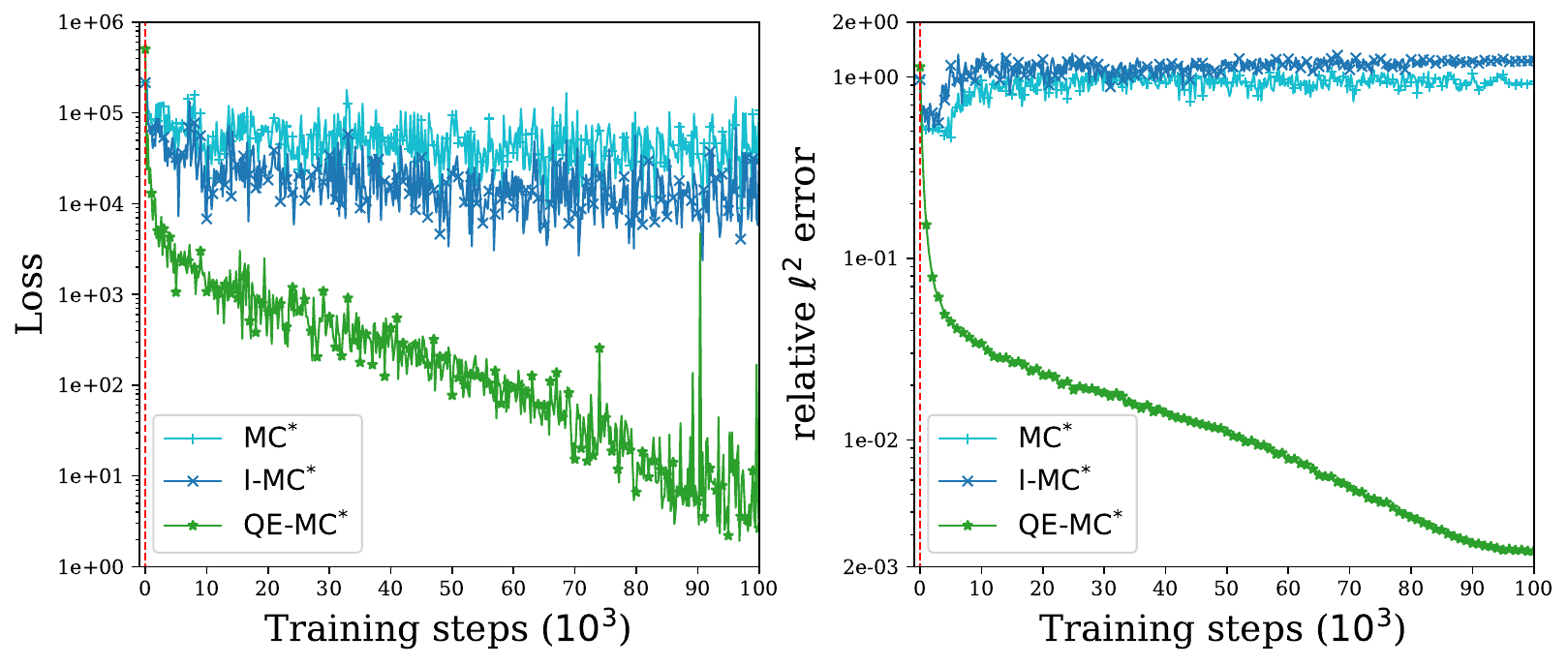}}
	
	\caption{Loss and relative $\ell^2$ error curves for the high-dimensional time-dependent fractional diffusion equation \cref{eq:fPDE} on the unit ball, with \(\{\gamma,\alpha,\boldsymbol{v},c\}=\{0.5,1.5,1,1\}\) and \(T=1\).}
	\label{fig:loss_error_fADE}
\end{figure}
\begin{table}[!htb]
	\centering
\caption{Relative $\ell^2$ errors $e_{\mathrm{test}}$ for the time-dependent fractional diffusion equation.}
	\label{tab:tfde_errors}
	\small
	\setlength{\tabcolsep}{5pt}
	\begin{tabular}{ccccccc}
		\toprule
		& \multicolumn{3}{c}{$u(\boldsymbol{x},t)$  in \cref{eq:time_sol_smooth}}
		& \multicolumn{3}{c}{$u(\boldsymbol{x},t)$  in \cref{eq:time_sol_singular}} \\
		\cmidrule(lr){2-4} \cmidrule(lr){5-7}
		$d$ & MC$^{*}$ & I-MC$^{*}$ & QE-MC$^{*}$ & MC$^{*}$ & I-MC$^{*}$ & QE-MC$^{*}$ \\
		\midrule
		$3$           &  3.01e-3   &  3.35e-2 & 1.77e-3  & 1.25e-1  & 4.88e-2  & 3.21e-3 \\
		$100$       &  9.07e-3  &  5.90e-2  & 9.38e-4 &  1.53e-1 & 9.93e-2 & 5.92e-3 \\
		$1{,}000$ &  1.70e-2   &  2.44e-1  & 4.84e-4 & 3.50e-1 & 5.21e-1   & 2.30e-3 \\
		$5{,}000$ &  2.92e-2  &  6.03e-1  & 2.82e-3 & 9.43e-1 & 1.23e00  & 2.44e-3  \\
		\bottomrule
	\end{tabular}
\end{table}

\cref{tab:tfde_errors,fig:loss_error_fADE} show that the proposed Quadrature-Enhanced MC-fPINN consistently provides the best accuracy and convergence behavior for the time-dependent fractional diffusion problem. For the smoother solution \cref{eq:time_sol_smooth}, Quadrature-Enhanced MC-fPINN attains the smallest relative \(L^2\) error for all tested dimensions. As the dimension increases, the errors of MC-fPINN and Improved MC-fPINN deteriorate much more rapidly than that of Quadrature-Enhanced MC-fPINN. For example, at \(d=1000\), the relative \(\ell^2\) errors of MC-fPINN, Improved MC-fPINN, and Quadrature-Enhanced MC-fPINN are \(1.70\times 10^{-2}\), \(2.44\times 10^{-1}\), and \(4.84\times 10^{-4}\), respectively; at \(d=5000\), the corresponding errors are \(2.92\times 10^{-2}\), \(6.03\times 10^{-1}\), and \(2.82\times 10^{-3}\). These results indicate that the proposed method maintains substantially better accuracy as the dimension grows.
For the more challenging solution \cref{eq:time_sol_singular}, which contains the lower-order boundary factor \((1-\|\boldsymbol{x}\|^2)^{\alpha/2}\), the advantage of the proposed method is even more pronounced. In this case, the relative \(\ell^2\) errors of MC-fPINN, Improved MC-fPINN, and Quadrature-Enhanced MC-fPINN are \(3.50\times 10^{-1}\), \(5.21\times 10^{-1}\), and \(2.30\times 10^{-3}\), respectively, at \(d=1000\); at \(d=5000\), they become \(9.43\times 10^{-1}\), \(1.23\times 10^{0}\), and \(2.44\times 10^{-3}\), respectively. Hence, while the two baseline methods suffer a substantial loss of accuracy as the dimension increases, Quadrature-Enhanced MC-fPINN still maintains errors at the level of \(10^{-3}\).

The convergence curves in \cref{fig:loss_error_fADE} are consistent with the error data in \cref{tab:tfde_errors}. In all cases, Quadrature-Enhanced MC-fPINN exhibits a more stable decay of both the loss and the relative \(\ell^2\) error, and converges to a significantly lower error level than MC-fPINN and Improved MC-fPINN. This advantage is particularly clear for the singular solution and becomes more significant in higher dimensions. Overall, these results demonstrate that the proposed quadrature-enhanced discretization substantially improves both robustness and accuracy for time-dependent fractional diffusion problems, especially in high dimensions and in the presence of strong boundary singularities.

\section{Conclusion}\label{sec:conclusion}
In this paper, we proposed a quadrature-enhanced Monte Carlo fPINN method for fractional Laplacian problems on bounded domains, including both the fractional Poisson equation and time-dependent fractional PDEs. The main difficulty of these problems lies in the simultaneous presence of hypersingular nonlocal operators, exterior Dirichlet constraints, reduced boundary regularity, and high dimensionality. To address these issues, we developed a geometry-adaptive discretization based on a spatially varying radius \(r_0(\boldsymbol{x})\) and directional distance-to-boundary information.

The proposed method decomposes the fractional Laplacian into three parts: a singular near-field term, a regular interior far-field term, and an exterior far-field term. In the numerical discretization, Gauss--Jacobi quadrature is employed for the singular radial integral, Gauss quadrature is used for the regular interior radial integral, and Monte Carlo sampling is adopted for the angular variables. This construction preserves the accuracy of deterministic radial quadrature while retaining the flexibility and scalability of Monte Carlo methods in high dimensions. To improve the approximation of low-regularity solutions, the discretization is further embedded into a feature-enhanced PINN trial space with explicit boundary factors and hard enforcement of the homogeneous initial and boundary conditions.

The numerical results demonstrate that the proposed method consistently outperforms the existing MC-fPINN and improved MC-fPINN discretizations in both accuracy and convergence behavior. In particular, the advantage is especially pronounced for bounded-domain problems with strong boundary singularities and for cases with relatively large fractional order. The method remains effective in very high dimensions; in particular, for the fractional Poisson equation, experiments are reported up to \(10{,}000\) dimensions, showing that the proposed quadrature-enhanced construction significantly improves robustness and approximation quality without sacrificing the mesh-free nature of the PINN framework.

Overall, the present work provides an accurate, robust, and scalable numerical framework for bounded-domain fractional Laplacian problems in high dimensions. Possible future directions include the extension to more general nonhomogeneous exterior conditions, tempered and anisotropic fractional operators, and more efficient training strategies for large-scale time-dependent problems.


\bibliographystyle{plain}
\bibliography{ref}

\end{document}